\documentclass[review]{elsarticle}

\usepackage[table,xcdraw]{xcolor}
\usepackage{graphicx}
\usepackage{amssymb}
\usepackage{amsmath} 
\usepackage{mathrsfs}
\usepackage{mathabx}
\usepackage{textcomp}
\usepackage{hyperref}
\usepackage{bbm}
\usepackage{pifont} 
\usepackage{color}
\usepackage{pgfplots}
\usepackage{amssymb}
\usepackage{amsmath}
\usepackage{multirow}
\usepackage{multicol}
\usepackage{subfig}
\usepackage{pifont}
\usetikzlibrary{datavisualization}
\usetikzlibrary{shapes.geometric}

\def\rojo{\color[rgb]{0.8,0.2,0.2}}

\def\IN{\mathbb N}

\def\1I{\mathbbm{1}}
\def\IZ{\mathbb Z}
\def\IR{\mathbb R}

\def\A{\mathcal A}

\def\D{\mathcal D}
\def\E{\mathcal E}

\def\N{\mathcal N}

\def\T{\mathcal T}

\def\e{\mathsf{e}}

\def\dist{{\rm dist}}

\def\ms{\medskip\noindent}

\journal{Physica D}

\begin{document}

\begin{frontmatter}

\title{Combinatorics of the paths towards synchronization}

\author[Fr,Mx]{A. Espa\~na}
\author[Fr]{X. Leoncini}
\author[Mx]{E. Ugalde}

\address[Fr]{Aix Marseille Univ, Université de Toulon, CNRS, CPT, Marseille, France.}
\address[Mx]{Instituto de F\'isica, Universidad Aut\'onoma de San Luis Potos\'i, M\'exico.}

\begin{abstract}
In this paper, we introduce a codification of the paths towards synchronization for synchronizing flows defined over a network. The collection of paths toward synchronization defines a combinatorial structure: the transition diagram. We describe the transition diagram corresponding to the Laplacian flow over the completely connected graph. This applies to the Kuramoto flow over the same graph when initial conditions close to the diagonal are considered. We present as well some results concerning the Laplacian and Kuramoto flows over the complete bipartite graph. 
\end{abstract}

\end{frontmatter}

%\linenumbers

\section{Introduction}

Synchronization phenomena are a long standing subject dating back at least to the observations of Huygens see for instance \cite{Pikovsky2001} This field of research when considering coupled dynamical systems on networks has been  very active since Kuramoto's seminal paper~\cite{Kuramoto1975}. The first studies considered homogeneously coupled systems like global coupling, completely random coupling or couplings according to a regular network. A very complete account of those early work can be found in~\cite{Strogatz2000}. As noticed in~\cite{Arenas&alPR2008}, the progression in connectivity of the synchronized subnetwork as times increases qualitatively follows the dictates of the linearized dynamics. Hence, the path towards synchronization can be understood through the study of the Laplacian as a linear dynamical system. The dynamics of the Laplacian is globally synchronizing and as we show below, the path towards the full synchronization can explicitly determined in that case. In contrast, the non-linear dynamics is not always fully synchronizing and some noticeable differences between linear and non-linear interactions appear as we increase the size of the system. Counting the number of paths to synchronization is a way of measuring the complexity of a system in the case of transient dynamics, and characterizing the complexity of the system by measuring the diversity of paths is a classic topic that has been studied and illustrated in \cite{Afraimovich2003,Afraimovich2005,Leoncini2002}.\\

The aim of this paper is to introduce the notion of synchronization sequences, which can be related to the connectivity matrix defined in~\cite{Arenas&al2006a}. In the case of a fully synchronizing system, the set of all the synchronizing sequences forms a transition diagram which encode the full transient dynamics towards synchronization. We study this combinatorial structure for the Laplacian dynamics on the complete graphs in full detail, and in some detail for the case of the complete bipartite graph. We characterize the transient dynamics on those networks by means of some topological features of the corresponding transition diagram.\\

The rest of the paper is organized as follows. After establishing the notations and the basic definitions in Section~\ref{sec:SetUp}, we study in Section~\ref{sec:KN} the transition diagram of synchronization paths for the complete graph $K_N$. Then, in Section~\ref{sec:KNN} we present some results concerning the structure of the transition diagram for the complete bipartite graph $K_{N,N}$. Finally, in Section~\ref{sec:Final} we close the paper with some final remarks and comments.

\section{Set-up}~\label{sec:SetUp}

We will refer indistinctly by graph or network to an undirected graph $G=(V,E)$, with vertices in $V$ and edges or links in $E$. 
On the contrary, a directed graph is a couple $D=(V,A)$ of vertices in $V$ and arrows in $A$. An edge is a two-vertex set, its end vertices, while an arrow is an ordered pair of vertices. A subgraph of $G$ is a graph $G'=(V',E')$ such that $V'\subset V$ and all the edges in $E'\subset E$ have end vertices in $V'$. A path in $G$ is a sequence of vertices such that each couple of consecutive vertices form and edge, while a path in $D$ is an ordered sequence of vertices $v_1\to v_2\to \cdots \to v_n$ such that each couple of consecutive vertices form an arrows. In this last case we say that $v_1$ is the starting vertex of the path and $v_n$ the ending one, besides  $n-1$, the number of arrows in the sequence, is the length of the path. A graph is connected if each couple of vertices belong to a path. Any graph can be decomposed in a unique way as a disjoint union of connected subgraphs, called connected components.\\

We fix a graph $G=(V,E)$ and consider a system of coupled differential equations on $I^V$, where $I$ is either $\IR$ or the circle $S^1$. The flow is generated by a system of ODEs coupled according to the edges in $E$ which represent the interactions between the particles. \\

We will focus on the Laplacian flow on $G$, which is the linear system defined by 
\begin{equation}
\frac{dx_v}{dt}=(L\,x)_v=\sum_{u\in V: \{u,v\}\in E}(x_u-x_v),\, \text{with } x_v\in \IR \,\text{ for each } v\in V.\label{eq:LaplacianFlow}
\end{equation}
Here $L$ is the Laplacian matrix of $G$, given by $L(v,v')=\sum_{u\in V: \{u,v\}\in E}(\1I_{\{u\}}-\1I_{\{v\}})(v')$. The synchronizing dynamics of the Laplacian flow is preserved in part by the Kuramoto flow defined in $(S^1)^V$ by the system of ODEs
\begin{equation}\label{eq:KuramotoModel}
\frac{dx_v}{dt}=\sigma\,\sum_{u\in V: \{u,v\}\in E}\sin (x_u-x_v),
\end{equation} 
where $\sigma\, \in\IR^+$ is the strength of the coupling. In both flows, the diagonal 
\begin{equation}\label{eq:Diagonal}
\D=\{x\in I^V:\, x_u=x_v\ \forall\ u,v\in V\},
\end{equation}
is an attractor, i.e., it is such that $\lim_{t\to\infty}\dist(x(t),\D)=0$, for each initial condition in a neighborhood of $\D$. Indeed, it is a global attractor for the Laplacian flow and, since the linearization of the Kuramoto flow around the diagonal is proportional to the Laplacian flow, applying a Hartman-Grobman argument we conclude that it follows a similar converging dynamics in a small neighborhood of the diagonal.  \\

\ms In order to measure the degree of synchronization at a given time, we fix a precision $\epsilon>0$ and declare that two neighboring sites are $\epsilon$-synchronized if their distance does not exceed $\epsilon$. Seeing as active each connection between neighboring sites which are $\epsilon$-close, we define a subnetworks containing all the active connections. The determination and evolution of this subnetwork is the main objective of the present work. 
Hence, to each fixed threshold $\epsilon > 0$ and every configuration $x\in\IR^V$, we associate an $\epsilon$-synchronized subnetwork $G_x=(V,E_x)$, where $E_x\subset E$ is the collection of edges 
\begin{equation}\label{eq:SubAdjacency}
E_{x}=\{\{u,v\}\in E:\, |x_u-x_v|\leq\epsilon\}.
\end{equation}
For the systems under consideration, $G_{x(t)}\rightarrow G$ as $t\to\infty$ provided the initial condition is sufficiently close to the diagonal. Since there is a finite number of subnetworks, then for each suitable initial condition $x\in\IR^V$ there exists a finite sequence of switching times $t_{0}=0 < t_{1}< t_{2} < \cdots < t_{\ell}$ and a corresponding sequence of $\epsilon$-synchronized subnetworks $(G_{x},G_{x(t_{1})},\ldots,G_{x(t_{\ell})})$ such that $G_{x(t_\tau)}\neq G_{x(t_{\tau+1})}$, for each $0\leq \tau < N$, and $G_{x(t)}=G_{x(t_{\tau})}$ with $\tau=\max\{0\leq j\leq \ell:\,t\geq t_{j}\}$. These sequence of subnetworks of $G$ codify the progression of transient synchronizing patterns. By taking $\epsilon$ sufficiently small, all the possible synchronizing sequences can be obtained by varying the initial condition $x$ inside the basin of attraction of the diagonal.

\ms In the case of highly symmetric networks, instead of the $\epsilon$-synchronized subnetworks it is convenient to use another combinatorial structure that at the same time that encodes the subnetwork and respects some of the symmetries that are preserved by the dynamics. As we will see below, this easy the description of the evolution of the synchronized subnetworks.
Hence, the whole synchronizing dynamics on $G$ can be compiled in a single combinatorial superstructure. This superstructure is a transition diagram whose vertices are in correspondence (not necessarily injective) with $\epsilon$-synchronized subnetworks in such a way that the collection of all the paths in the transition diagram is equivalent to the set of all the observable sequences of $\epsilon$-synchronized subnetworks. To study this dynamic, it is enough to see the diagram with other labels that allow to encode the $G_x$. To be more precise, the transition diagram is a directed graph $\T_\epsilon=(V_\epsilon,A_\epsilon)$ whose vertices, $V_\epsilon$ are combinatorial objects containing all the information we need to determine the $\epsilon$-synchronized subnetworks and the arrows, $A_\epsilon$, are transitions between those structures which are consistent with the evolution of the synchronized subnetworks. The correspondence between objects in $V_\epsilon$ and $\epsilon$-synchronized subnetworks is achieved via a mapping
\begin{equation}\label{eq:transitiondiagram}
\lambda:V_\epsilon\to \E_\epsilon,
\end{equation} 
that labels each vertex in the transition diagram with an $\epsilon$-synchronized subnetwork. The labelling $\lambda$ is such that $(G_0,G_1,\ldots,G_\ell)$ is a realizable sequence of $\epsilon$-synchronized subnetworks if and only if there exists a path $v_0\to v_1\to \cdots \to v_\ell$ in $\T_\epsilon$ such that $G_n=\lambda(v_n)$ for $0\leq n\leq \ell$, we will see these encodings in detail later. \\

In general, the set $\E_{\epsilon}$ of all the $\epsilon$-synchronized subnetworks changes with $\epsilon$. Nevertheless, for $\epsilon$ sufficiently small, the set of $\epsilon$-synchronized subgraphs defined by initial conditions in a small neighborhood of $\D$ becomes independent of $\epsilon$. For the Laplacian flow, the set $\E_{\epsilon}$ of all possible synchronized subnetworks is independent of $\epsilon$ as long as $\epsilon>0$. Even if $\E_\epsilon$ is independent of $\epsilon$, the corresponding transition diagram may change with $\epsilon$. This, nevertheless, does not happen in the linear case, since for each initial condition $x\in \IR^V$, the corresponding sequence $(G_{x},G_{x(t_1)},\ldots,G_{x(t_{\ell})})$ of $\epsilon$-synchronized subnetworks coincides with the sequence $(G_{y},G_{y(t_1)},\ldots,G_{y(t_{\ell})})$ of $\epsilon'$-synchronized subnetworks determined by $y=x\,\epsilon'/\epsilon$.  Indeed, by Equation~(\ref{eq:SubAdjacency}) and by the linearity of the system, $\{u,v\}\in E_x$ is equivalent to $|x_u-x_v|\leq \epsilon$, hence $|x_u-x_v|=\epsilon/\epsilon'|y_u-y_v|\leq \epsilon$, therefore $|y_u-y_v|\leq \epsilon'$, which is equivalent to  $\{u,v\}\in E_y$. From this it follows that the collection of $\epsilon$-synchronized sequences does not depend on $\epsilon$ in the linear case. Clearly, since this number is finite, each synchronized  sequence can be realized by an infinite number of initial conditions, which could most likely allow to realize some partition of the initial phase space, i.e.,  the basin of attraction of the final synchronized state.\\

We will restrict our study to the following families of networks:
\begin{itemize}
\item[A.] The complete graph $K_N$, for which $V=\{1,2,\ldots,N\}$ and $E=\{\{u,v\}: 1\leq u < v\leq N\}$.
\item[B.] The complete bipartite graph $K_{N,N}$, where $V=\{1,2,\ldots,2N\}$ and $E=\{\{u,N+v\}:\, 1\leq u,v\leq N\}$. 
\end{itemize}
Considering these families,  we address the following questions:
\begin{enumerate}
\item Given the underlying network, which subgraphs are realizable as synchronizing subnetworks? How large is this collection and how does it grow with the size of the underlying graph?
\item Given an underlying network, what is the structure of the transition diagram? In particular, what is the longest path in this digraph and what is the resulting distribution of path lengths? 
\end{enumerate}

\ms 

\section{The transition diagram for $K_N$}\label{sec:KN}
\ms 
The Laplacian matrix for $K_N$ diagonalizes in the basis $\{u^1,u^2,\ldots,u^N\}$, where $u^1:=\sum_{n=1}^N\e^n$ and, for each $n\geq 1$, $u^n:=\e^n-\e^1$, where $\e^n$ denotes the $n$-th vector of the canonical basis. Indeed, $Lu^1=0$ and $Lu^n=-N\,u^n$ for each $n\geq 2$. Consider now an initial condition $x\in \IR^N$. Such an initial condition can be decomposed as $x=\bar{x}\,u^1+\sum_{n=1}^{N-1}(x_{n+1}-\bar{x})\,u^n$, where $\bar{x}:=\left(\sum_{n=1}^Nx_n(0)\right)/N$. Therefore
\[
x(t)=\bar{x}\,u^1+ e^{-N\,t}\sum_{n=1}^{N-1}(x_{n+1}-\bar{x})\,u^n
=\sum_{n=1}^N\left(\bar{x}
\left(1-e^{-N\,t}\right)+e^{-N\,t}x_n\right)\e^n,
\]
for all $t\in\IR$. From this it follows that
\begin{equation}~\label{eq:KNMonotonicity}
x_n(t)-x_m(t)=e^{-N\,t}(x_n-x_m), 
\end{equation}
for all $t\in\IR$ and each $1\leq m,n\leq N$. Hence, the edge $\{n,m\}$ belong to the $\epsilon$-synchronized subnetwork $G_{x(t)}$, for all times exceeding 
$t_{n,m} = \left(\log|x_n-x_m|-\log(\epsilon)\right)/N$.

\ms Without lost of generality we may assume that $x_1(0)\leq x_2(0)\leq\cdots\leq x_N(0)$ which, by Equation~\eqref{eq:KNMonotonicity}, ensures that $x_1(t)\leq x_2(t)\leq\cdots\leq x_N(t)$ for all $t$.

\ms In order to take advantage of the fact that the Laplacian flow preserves the order of the coordinates, we will define the transition diagram not over the synchronized subnetworks but over another combinatorial object that encodes both the synchronized subnetworks, and recognizes the order of the coordinates. By doing so we will facilitate the description of the transition diagrams since the coding we use allows us to easily determine the order of apparition of new edges in the synchronized sequence. This coding is not only convenient but necessary if one wants to keep track of the order of the coordinates. We codify $\epsilon$-synchronized subnetwork $G_x$, determined by the ordered configuration $x_1\leq x_2\leq\cdots\leq x_N$ by the increasing function $\phi_x:\{1,2,\ldots,N\}\to\{1,2,\ldots,N\}$ given by
\begin{equation}\label{eq:KNIncreasingFunction}
\phi_x(m)=\max\{n\geq m:\, x_n\leq x_m+\epsilon\}.
\end{equation}
Clearly $\phi_x$ is increasing and such that $\phi_x(n)\geq n$ for each $1\leq n\leq N$, i.e., $\phi_x\geq {\rm Id}$. Here and below Id denotes the identity function in $\{1,2,\ldots,N\}$. We present an example of the construction of the increasing function from a given initial condition, in Figure~\ref{fig:KN4ConstructionEx}.

\ms By the arguments in the Appendix~\ref{app:KNIncreasing}, the collection 
\begin{equation}\label{eq:KNFunctionsCollection}
\Phi_N:=\{\phi:\{1,\ldots,N\}\to\{1,\ldots,N\}\text{ increasing and such that } \phi\geq {\rm Id}\},
\end{equation}
is in a one-to-one correspondence with the collection of all $\epsilon$-synchronized subnetworks of $K_N$ defined by initial conditions satisfying $x_1\leq x_2\leq\cdots\leq x_N$. The correspondence is given by
\begin{equation}~\label{eq:KNNetworkFuctionCorrespondance}
\phi\mapsto (\{1,2,\ldots,N\}, E_\phi)\, \text{ where }
E_{\phi}=\{\{m,n\}:\,\min(m,n)\leq \phi(\max(n,m))\}.
\end{equation}

In this case, the coding~\eqref{eq:transitiondiagram} which associates increasing functions to synchronized subnetworks is given by Equation~(\ref{eq:KNNetworkFuctionCorrespondance}).

\begin{center}
\begin{figure}[h]
\begin{tikzpicture}

%%%%%%%%%% Vecindades
\draw  (0,0)  -- (10,0);
\draw  (1,0.5)--(3.5,0.5)--(3.5,-1)--(1,-1)--(1,0.5);
\draw  (6.25,0.5)--(8.75,0.5)--(8.75,-1)--(6.25,-1)--(6.25,0.5);

%%%%%%%%%% Coordinates
\filldraw (8,0)    circle (3pt);
\filldraw (7,0)    circle (3pt);
\filldraw (3,0)    circle (3pt);
\filldraw (1.5,0)  circle (3pt);

%%%%%%%%%% Coordinates labels
\draw (8,-0.5) node {\text{$x_4$}};
\draw (7,-0.5) node {\text{$x_3$}};
\draw (3,-0.5) node {\text{$x_2$}};
\draw (1.5,-0.5) node {\text{$x_1$}};

%%%%%%%%%% Arrows
\draw (2.25,-1.5) node {\text{$\downarrow$}};
\draw (5,-1.5) node {\text{$\downarrow$}};
\draw (7.5,-1.5) node {\text{$\downarrow$}};

%%%%%%%%%% Arrow labels
\draw (2.25,-2) node {\text{Link}};
\draw (5,-2) node {\text{No link}};
\draw (7.5,-2) node {\text{Link}};

%%%%%%%%%% Vertices
\filldraw [color=blue] (9.5,-3) circle (3pt);
\filldraw [color=blue] (6.5,-3)  circle (3pt);
\filldraw [color=blue] (3.5,-3)  circle (3pt);
\filldraw [color=blue] (0.5,-3)  circle (3pt);

%%%%%%%%%% Vertex labels
\draw (9.5,-3.5) node {\text{4}};
\draw (6.5,-3.5) node {\text{3}};
\draw (3.5,-3.5) node {\text{2}};
\draw (0.5,-3.5) node {\text{1}};

%%%%%%%%%% Edges
\draw (9.5,-3) -- (6.5,-3);
\draw (0.5,-3) -- (3.5,-3);

%%%%%%%%%% Text Node
\draw (5,-5) node {\text{$\phi_x$=(2,2,4,4)}};

%%%%%%%%%% Label Text 
\draw (-1,0)  node {\text{(a)}};
\draw (-1,-3) node {\text{(b)}};
\draw (-1,-5) node {\text{(c)}};
\end{tikzpicture}
\caption{In (a), an example of the values of $x=(x_1,x_2,x_3,x_4)$ are illustrated with black dots. To construct $G_{x}$, according to Equation~(\ref{eq:SubAdjacency}), it is enough to observe that $x_1$ and $x_2$ are inside one $\epsilon$-neighborhood, and $x_3$ and $x_4$ in another, which implies that in (b) there are a links between the vertices 1 and 2 as well as vertices 3 and 4. In (c), the increasing function $\phi_x$ associated with $x$ is depicted. The information in $\phi_x$ can be read as follows: The furthest vertex connected with vertex 1 is vertex 2, vertex 2 does not reach vertex 3, and vertex 3 reaches vertex 4, which is the last one.}~\label{fig:KN4ConstructionEx}
\end{figure}
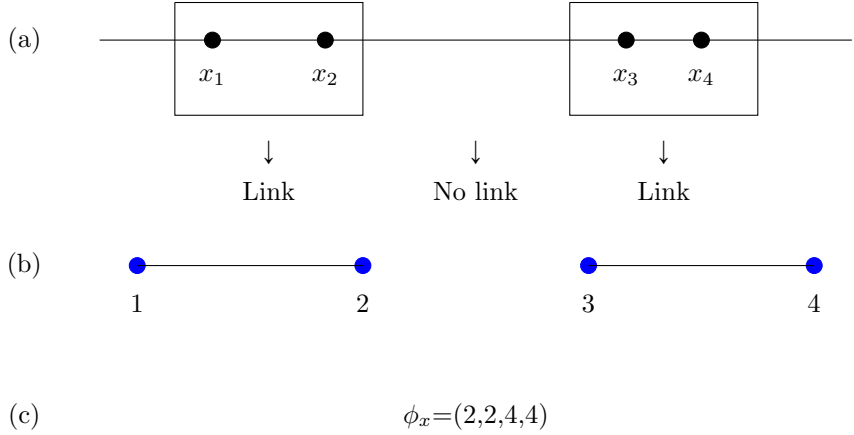
\end{center}

The collection $\Phi_N$ is equivalent to a well-studied combinatorial set, the set $\D_{2N}$ of Dyck paths of length $2N$. This set is in turn equivalent to the set of valid $2N$-parenthesis configurations. All these combinatorial sets have a cardinality given by the Catalan numbers~\cite{Stanley1999}, 
\begin{equation}\label{eq:Catalan}
|\Phi_N|=C_N:=\frac{1}{N+1}
\left(\begin{matrix}2\,N\\N\end{matrix}\right). 
\end{equation} 

\ms Taking into account the equivalence established in the previous paragraph, each sequence of $\epsilon$-synchronized subnetworks $(G_{x},G_{x(t_1)},\ldots,G_{x(t_\ell)})$ generated by an ordered initial condition $x$, is faithfully codified by the corresponding sequences of increasing functions $(\phi_{x},\phi_{x(t_1)},\ldots,\phi_{x(t_\ell)})$ in $\Phi_N$ defined by Equation~\eqref{eq:KNIncreasingFunction}. Clearly the function $t\mapsto\phi_{x(t)}(n)$ increases with $t$ for each $n$ fixed, and converges to $\phi_{x(t)}=N$ at the time $t_{1,N} = \left(\log(x_N-x_1)-\log(\epsilon)\right)/N$. Due to the monotonicity, the length $\ell$ of an $\epsilon$-synchronized sequences is upper bounded by the number of edges in $K_{N}$, i.e., $\ell\leq N(N-1)/2$. As mentioned above, the switching times $t_1<t_2<\cdots<t_{\ell}$ are completely determined by the increments $x_n-x_m$, with $m < n$. Let us assume that all those increments are different from zero and pairwise different. We will say that a path satisfying this condition is typical. Clearly, the non-typical paths correspond to initial conditions in a set of zero Lebesgue measure in $\IR^N$. Hence, for typical paths, two consecutive functions in the sequence $(\phi_0,\phi_1,\ldots,\phi_\ell):=(\phi_{x},\phi_{x(t_1)},\ldots,\phi_{x(t_\ell)})$ differ at a single point. Let us denote by $\delta_n\in\{0,1\}$ the characteristic function of the singleton $\{n\}$. Hence $\phi_{\tau+1}=\phi_\tau+\delta_{n_\tau}$ for some $n_\tau\in\{1,2,\ldots,N\}$ satisfying the condition $\phi_\tau(n_\tau)<\phi_\tau(n_\tau+1)$. Hence an admissible sequence $(\phi_0,\phi_1,\ldots,\phi_\ell)$ can be obtained by choosing a valid initial function $\phi_0\in\Phi_N$, then for each $\tau\geq 0$, a point $n_\tau\in\{1,2,\ldots,N-1\}$ such that $\phi_\tau(n_\tau)<\phi_\tau(n_\tau+1)$ to update $\phi_{\tau+1}=\phi_\tau+\delta_{n_\tau}$. Nevertheless, not all the sequences obtained in this way are realizable as synchronizing sequences. The sequence $\left(n_\tau\right)_{0\leq\tau<\ell}$ of jump sites is determined by an order in the increments $\Delta:=\{\Delta_{n,k}:=x_{n+k}-x_n:\, 1\leq n < n+k\leq N\}$ in such a way that the $\tau$-th smallest increment in $\Delta$ is of the kind $\Delta_{n_{\tau},k}$. Hence, to each valid strict ordering in $\Delta$ corresponds a unique realizable path towards synchronization. 

\ms One can easily verify that not all the admissible paths are realizable. The simplest counterexample happens for $N=4$ (for $N=2,3$ all admissible sequences are realizable). In this case the sequence ${\rm Id}\mapsto (2,2,3,4) \mapsto (2,2,4,4)\mapsto (2,3,4,4)\mapsto (2,4,4,4)\mapsto (3,4,4,4)\mapsto (4,4,4,4)$, which corresponds to the sequences of jump sites $(1,3,2,2,1,1)$, is not realizable since the first two transitions indicate that $x_2-x_1 < x_4-x_3$ but transitions four and five imply that $x_4-x_2 < x_3-x_1$, which is contradictory. 
The total number of admissible paths  for $N=4$ is sixteen. On the other hand, the total number of admissible paths is ten, and the associated valid strict orderings are shown in Table~\ref{tab:KN4Ordering}.

\begin{table}[h]
\begin{tabular}{|c|c|}
\hline
$\Delta_{1,1} < \Delta_{2,1} < \Delta_{3,1} < \Delta_{1,2} < \Delta_{2,2} < \Delta_{1,3}$ &
$\Delta_{1,1} < \Delta_{2,1} < \Delta_{1,2} < \Delta_{3,1} < \Delta_{2,2} < \Delta_{1,3}$ \\
$\Delta_{1,1} < \Delta_{3,1} < \Delta_{2,1} < \Delta_{1,2} < \Delta_{2,2} < \Delta_{1,3}$ &
$\Delta_{2,1} < \Delta_{1,1} < \Delta_{3,1} < \Delta_{1,2} < \Delta_{2,2} < \Delta_{1,3}$ \\
$\Delta_{2,1} < \Delta_{1,1} < \Delta_{1,2} < \Delta_{3,1} < \Delta_{2,2} < \Delta_{1,3}$ &  
$\Delta_{2,1} < \Delta_{3,1} < \Delta_{1,1} < \Delta_{2,2} < \Delta_{1,2} < \Delta_{1,3}$ \\
$\Delta_{2,1} < \Delta_{3,1} < \Delta_{2,2} < \Delta_{1,1} < \Delta_{1,2} < \Delta_{1,3}$ &
$\Delta_{3,1} < \Delta_{1,1} < \Delta_{2,1} < \Delta_{2,2} < \Delta_{1,2} < \Delta_{1,3}$ \\
$\Delta_{3,1} < \Delta_{2,1} < \Delta_{1,1} < \Delta_{2,2} < \Delta_{1,2} < \Delta_{1,3}$ &
$\Delta_{3,1} < \Delta_{2,1} < \Delta_{2,2} < \Delta_{1,1} < \Delta_{1,2} < \Delta_{1,3}$\\
\hline
\end{tabular}
\ms 
\caption{The ten different orderings of the increments for a typical initial conditions in $\IR^4$.}\label{tab:KN4Ordering}
\end{table}

\ms Each ordering in Table~\ref{tab:KN4Ordering} uniquely determines an observable path towards synchronization. The corresponding paths towards synchronization are organized in a transition diagram, depicted in Figure~\ref{fig:KN4TransitionDiagram}.

\begin{center}
\begin{figure}[h]
\tikzset{every picture/.style={line width=0.7pt}}     

\begin{tikzpicture}[x=0.75pt,y=0.75pt,yscale=-1,xscale=1]

%%%%%%%%%% Arrows
\draw [-to] (367,25)  -- (285,98)  ;
\draw [-to] (367,25)  -- (444,98)  ;
\draw [-to] (367,25)  -- (367,98)  ; 
\draw [-to] (285,104) -- (285,180) ;
\draw [-to] (448,104) -- (448,180) ;
\draw [-to] (367,106) -- (441,180) ;
\draw [-to] (367,104) -- (367,180) ;
\draw [-to] (285,104) -- (362,180) ;
\draw [-to] (448,104) -- (525,180) ;
\draw [-to] (285,187) -- (285,262) ;
\draw [-to] (367,187) -- (289,262) ;
\draw [-to] (367,187) -- (367,262) ;
\draw [-to] (448,187) -- (448,262) ;
\draw [-to] (527,187) -- (527,262) ;
\draw [-to] (448,187) -- (522,262) ;
\draw [-to] (367,270) -- (367,348) ;
\draw [-to] (285,270) -- (362,348) ;
\draw [-to] (448,270) -- (448,348) ;
\draw [-to] (527,270) -- (452,348) ;
\draw [-to] (448,355) -- (448,430) ;
\draw [-to] (366,355) -- (442,430) ;
\draw [-to] (448,437) -- (448,512) ;
%%%%%%%%%%%%%%%%%%%%%%%%%%%%%%%%%%%%%%%%%%%%%%%

%%%%%%%%%%%%% Vertices
\filldraw [color=blue] (285,104) circle (3pt);
\filldraw [color=blue] (285,187) circle (3pt);
\filldraw [color=blue] (285,270) circle (3pt);
\filldraw [color=blue] (367,23)  circle (3pt);
\filldraw [color=blue] (367,104) circle (3pt);
\filldraw [color=blue] (367,187) circle (3pt);
\filldraw [color=blue] (367,270) circle (3pt);
\filldraw [color=blue] (367,355) circle (3pt);
\filldraw [color=blue] (448,104) circle (3pt);
\filldraw [color=blue] (448,187) circle (3pt);
\filldraw [color=blue] (448,270) circle (3pt);
\filldraw [color=blue] (448,354) circle (3pt);
\filldraw [color=blue] (448,437) circle (3pt);
\filldraw [color=blue] (448,520) circle (3pt);
\filldraw [color=blue] (528,185) circle (3pt);
\filldraw [color=blue] (528,270) circle (3pt);
%%%%%%%%%%%%%%%%%%%%%%%%%%%%%%%%%%%%%%%%%%%%%%%%%%%%%

%%%%%%%%%%%%%% Vertex labels
\draw (382,20)  node {\text{Id}};
\draw (322,100) node {\text{(2,2,3,4)}};
\draw (402,100) node {\text{(1,3,3,4)}};
\draw (487,100) node {\text{(1,2,4,4)}};
\draw (322,183) node {\text{(2,2,4,4)}};
\draw (402,183) node {\text{(2,3,3,4)}};
\draw (487,183) node {\text{(1,3,4,4)}};
\draw (561,183) node {\text{(2,2,4,4)}};
\draw (561,265) node {\text{(2,3,4,4)}};
\draw (487,265) node {\text{(1,4,4,4)}};
\draw (402,265) node {\text{(3,3,3,4)}};
\draw (322,265) node {\text{(2,3,4,4)}};
\draw (402,351) node {\text{(3,3,4,4)}};
\draw (487,351) node {\text{(2,4,4,4)}};
\draw (487,433) node {\text{(3,4,4,4)}};
\draw (463,516) node {\textbf{4}};
%%%%%%%%%%%%%%%%%%%%%%%%%%%%%%%%%%%%%%%%%%%%%
\end{tikzpicture}
\caption{The transition diagram which contains all the paths towards synchronization of the Laplacian dynamics on $K_4$. The synchronized subgraphs are encoded by increasing functions as defined by Equation~\eqref{eq:KNNetworkFuctionCorrespondance}. At the top is placed the identity function ${\rm Id}:=(1,2,3,4)$ which codifies the completely disconnected graph. All the paths end at the constant function ${\bf 4}=(4,4,4,4)$,  which codifies the globally synchronized state}.~\label{fig:KN4TransitionDiagram}
\end{figure}
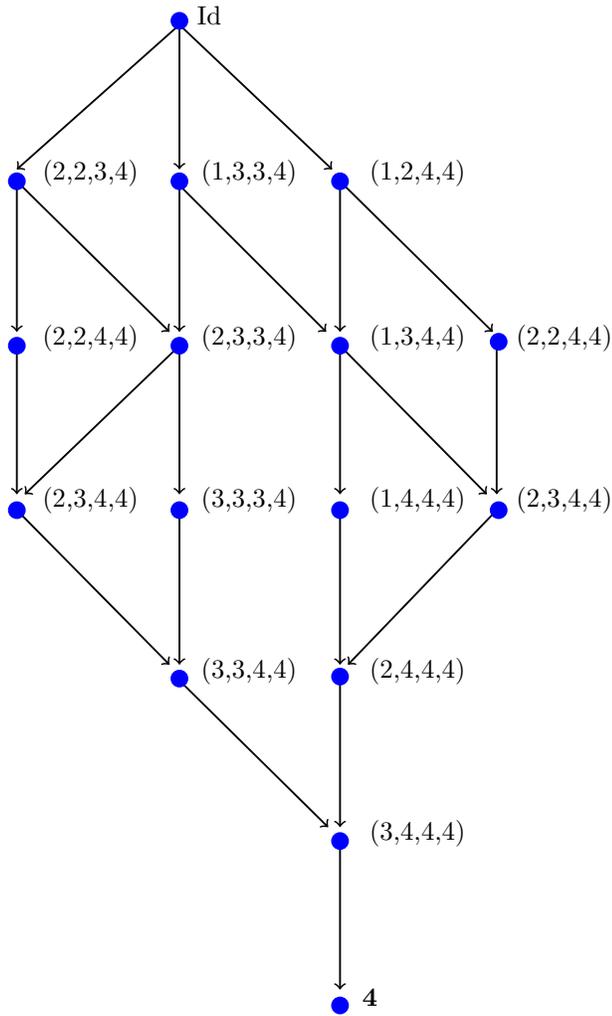 
\end{center}

\ms As mentioned above, the path towards synchronization form the initial condition $x$ is given by the sequence $(G_{x},G_{x(t_1)},\ldots, G_{x(t_\ell)})$ of synchronizing subnetwork, which is equivalent to a sequences of increasing functions $(\phi_{x},\phi_{x(t_1)},\ldots, \phi_{x(t_\ell)})$ in $\Phi_N$. The sequence $(\phi_{x},\phi_{x(t_1)},\ldots, \phi_{x(t_\ell)})$ is completely determined by the order of the increments $\Delta$.  Each ordering of increments determines the sequence $\left(n_\tau\right)_{0\leq\tau<\ell}$ of sites where consecutive increasing functions differ, i.e., the sites $n_\tau$ such that $\phi_{x(t_{\tau+1})}-\phi_{x(t_{\tau})}=\delta_{n_\tau}$ for each $0\leq \tau < \ell$. Hence, each valid ordering in $\Delta$ corresponds a unique realizable path towards synchronization. Therefore, the total number of paths toward synchronization is given by the number of different orderings $\Delta$ which can be obtained from an ordered vector $x\in\IR^N$. This is a combinatorial problem which has been treated in the literature in the context of the so called Golomb rulers \cite{Golomb1972}, that is, the problem of counting the number of valid orders is equivalent to counting the combinatorially distinct Golomb rulers. Below we will explain how this equivalence is established.

\ms A Golomb ruler with $N$ marks is a vector $a\in \IZ^N$ with $a_1 < a_2 < \cdots < a_N$, such that no two increments  $a_{n+k}-a_{n}$, where $1\leq n< N$, and $1\leq k\leq N-n$ coincide. Hence, a Golomb ruler is nothing but a typical initial condition with integer entries. 

\ms To each typical initial condition $x\in\IR^N$ we may associate a Golomb ruler as follows. Since $x$ is typical, then both
$\epsilon_1=\min\{\Delta_{n,k}\,:\, 1\leq n < N,\, 1\leq k <N-n\}$ and $\epsilon_2=\min\{|\Delta_{n,k}-\Delta_{m,\ell}|\,:\, (m,k)\neq (n,\ell) \,:\, 1\leq n < N,\, 1\leq k <N-n , \, 1\leq m < N,\, 1\leq \ell <N-m \}$ are strictly positive. Let $p\in\IN$  be such that $p\cdot\min(\epsilon_1,\epsilon_2/4) > 1$, and for each $1\leq n \leq N$ let $q_n:=\max\{q\in\IZ:\, q/p \leq x_n\}$. The vector $q=(q_1,q_2,\ldots,q_N)\in\IZ^N$ is the desired Golomb ruler. Indeed, since $p\,\epsilon_1 > 1$, then for each $1\leq n <N$ we have 
\[q_n\leq p\,x_n \leq p\,(x_{n+1}-\epsilon_1)\leq q_{n+1}+1-p\,\epsilon_1 < q_{n+1}.\] 
On the other hand, whenever $\Delta_{n,k} > \Delta_{m,\ell}$ we have 
\[(q_{n+k}-q_n)-(q_{m+\ell}-q_m) \geq p(\Delta_{n,k}-\Delta_{m,\ell}-4/p) > p(\epsilon_2 - 4/p) > 0.\]

\ms Two Golomb rulers are combinatorially equivalent if they determine the same ordering in their differences, i.e., $a,b\in\IR^N$ are equivalent if and only if $((a_{n+k}-a_n)-(a_{m+\ell}-a_m))((b_{n+k}-b_n)-(b_{m+\ell}-b_m)) >0$ for each $1\leq m,n < N$ and $1\leq k < n,\, 1\leq \ell < m$. Hence, the number ${\rm Golomb}(N)$ of classes of Golomb rulers with $N$ marks, gives the number of paths towards synchronization, i.e.,
\begin{equation}\label{eq:NoPathsKN}
\text{Number of paths towards synchronization for }K_N=\text{Golomb}(N).
\end{equation} 

The growth of this quantity with the dimension $N$, is a measure of complexity similar to the topological complexity of discrete-time dynamical systems. In the case of a discrete-time dynamical system, the topological complexity counts the growth of the number of distinguishable trajectories as a function of time. In our case,  $\text{Golomb}(N)$ counts the number of distinguishable paths towards synchronization, not as a function of time, but of the dimension of the system.

\ms A Golomb ruler $a\in\IZ$ is also characterized by the fact that all the sums $a_m+a_n$ are different. Indeed, since
\[
{\rm sign}((a_{n+k}-a_n)-(a_{m+\ell}-a_m))
={\rm sign}((a_{n+k}+a_m)-(a_{m+\ell}+a_n)),
\] 
the number of combinatorially different Golomb rules is given by the number of different orderings for $S=\{a_m+a_n:\, 1\leq m < n\leq N\}$ which is equal to the number of different orderings for $P=\{a_m\,a_n:\, 1\leq m < n\leq N\}$. This number is relevant in problem of quantum entanglement~\cite{Hildebrand2002}. The sequence ${\rm Golomb}(N)$ appears in the On-line Encyclopedia of Integer Sequences under the entry A237749~\cite{web:oeisA237749}, where the first nine terms, which we present in Table~\ref{tab:Golomb}, are explicitly computed. 

\begin{table}[h!]
\centering
\begin{tabular}{|c|l|}
\hline
$N$ & Golomb$(N)$ \\ 
\hline
1          & 1         \\ \hline
2          & 1         \\ \hline
3          & 2         \\ \hline
4          & 10        \\ \hline
5          & 114       \\ \hline
6          & 2608      \\ \hline
7          & 107498    \\ \hline
8          & 7325650   \\ \hline
9          & 771505180 \\ \hline
\end{tabular}

\ms
\caption{Number of classes of Golomb rulers.}\label{tab:Golomb}
\end{table}

\ms The computation of ${\rm Golomb}(N)$ remains an open problem.  Easy bounds for this number are shown in Equation~(\ref{eq:BoundGolomb}). The lower bound can be obtained by counting all the orderings of the first differences $x_{i+1}-x_i$ for $1\leq i\leq N-1$, while the upper bound results taking all the ordering of all the differences $x_{i} -x_k$ for $1\leq k<i\leq N$. Form this we obtain,
\begin{equation}~\label{eq:BoundGolomb}
(N-1)!< {\rm Golomb}(N) <\binom{N}{2}!.
\end{equation}

\ms An exact non-trivial upper bound, based on a result by M. R. Thrall~\cite{Thrall1952}, was found by N. Johnston~\cite{webJohnston}. It establishes that
\begin{equation}~\label{eq:UpperBoundGolomb}
{\rm Golomb}(N) \leq \frac{\prod_{n=1}^{N-1}n!}{\prod_{n=1}^N(2n-1)!}\  \left( \frac{N(N+1)}{2}\right)!\,
\end{equation}
which furnishes an upper bound for the number of paths towards synchronization as well. 

\ms We do not intend to make an exhaustive characterization of the transition diagram, but from the concepts already defined, certain characteristics can be calculated, such as: the number of synchronized sequences of length $\ell$, the distribution of lengths of the  path towards synchronization. From this we compute the mean length and the most frequent length of paths. Furthermore, we can extrapolate the behavior of these quantities for increasingly large dimensions.

\ms The transition diagram for $K_N$ has a hierarchical structure with the  disconnected subnetwork, codified by the identity function ${\rm Id}\in\Phi_N$, at the top, and the completely connected network, codified by the constant function ${\bf N}(n)=N$, at the bottom. Since we are considering only typical initial conditions, at each transition only one new edge appears in the $\epsilon$-synchronized subnetwork. At level $\ell$, from top to bottom, we place all the subnetworks which can be reached from the disconnected subnetwork after exactly $\ell$ transition. These subnetworks are precisely those having exactly $n$ edges, and are therefore codified by increasing functions $\phi\in\Phi_{N}$ such that $\sum_{n=1}^N(\phi(n)-n)=\ell$. In particular, the maximal length of a synchronizing sequences is $l_{\max}=\sum_{n=1}^N(N-n)=N(N-1)/2$. From our discussion above is readily follows that the number $F_N(\ell)$ of synchronized sequences of length $\ell$ is given by number of Dyck paths of length $2N$ and area $N^2-\ell$, i.e.,
\begin{equation}\label{eq:DistLengths}
F_N(\ell):=\left|\left\{\phi\in\Phi_N:\ \sum_{n=1}^N\phi(n)=N^2-\ell\right\}\right|.
\end{equation}
These quantities can be computed from the generating polynomials 
\[P_N(t):=\sum_{\phi\in\Phi_N}t^{{\rm area}(\phi)}=
\sum_{\ell=0}^{\frac{N(N-1)}{2}} F_N(\ell)\,t^{\frac{N(N-1)}{2}-\ell},\] 
where ${\rm area}(\phi)=\sum_{n=1}^N(\phi(n)-n)$ denotes the area under the Dyck path determined by the increasing function $\phi$. The generating polynomials can be determined by using the recurrence relation 
\begin{equation}~\label{eq:RecurrenceCarliz}
P_N(t)=\sum_{n=0}^{N-1}t^n\,P_n(t)\,P_{N-n-1}(t)
\end{equation}
with initial conditions $P_0=0$, derived by Carlitz and Riordan~\cite{Carliz&Riordan1964} (see~\cite{BlancoPetersen2014} as well). Although there is no closed form for $F_N(\ell)$, the recurrence relation above allows to directly compute these distributions and to establish its asymptotic behavior. In Table~\ref{tab:Carliz} we show $F_N(\ell)$ for $2\leq n\leq 8$.

\ms
\begin{table}[h!]
\centering
\begin{tabular}{|c|l|}
\hline
$N$ & $F_N(\ell)$                  \\ \hline
2 & (1,1)                          \\ \hline
3 & (1,1,2,1)                      \\ \hline
4 & (1,1,2,3,3,3,1)                \\ \hline
5 & (1,1,2,3,5,5,7,7,6,4,1)        \\ \hline
6 & (1,1,2,3,5,7,9,11,14,16,16,17,
14,10,5,1)                         \\ \hline
7 & (1,1,2,3,5,7,11,13,18,22,28,32,37,40,44,43,40,35,25,15,6,1)                      \\ \hline
8 & (1,1,2,3,5,7,11,15,20,26,34,42,53,63,73,85,96,106,113,118,118,115,102,86,65,41,21,7,1)\\ \hline
\end{tabular}

\ms
\caption{Number $F_{N}(\ell)$ of functions $\phi\in\Phi_{N}$ codifying a subnetworks starting a synchronizing path of length $\ell$.}~\label{tab:Carliz}
\end{table}

\ms The normalized cumulative distribution, 
$f_N:[0,1]\to[0,1]$, is defined by 
\begin{equation}~\label{eq:ComulativeDist}
f_N(x)=\frac{1}{C_N}
\sum_{n\leq x\times N(N-1)/2}F_N(x),
\end{equation}
where $F_N$ is given by Equation~\eqref{eq:DistLengths} and $C_N$ the $N$-th Catalan number. By using the recurrence shown in Equation~\eqref{eq:RecurrenceCarliz}, we numerically computed $f_N(x)$ for increasing values of $N$, and observe that $f_N$ approaches an absolutely continuous limit distribution $x\mapsto f(x)$ whose density $\rho(x):=d\,f(x)/dx$ is closely approached by the curve depicted in Figure~\ref{fig:DistributionLength}. Hence, for $N$ sufficiently large and $\delta>0$ sufficiently small, the proportion of paths of length $N(N-1)(x\pm \delta)/2$ is approximatively $\rho(x)\,\delta$. As shown in the figure, our numerical computation suggest that $\rho$ is continuous, unimodal, and negatively skewed. 

\begin{center}
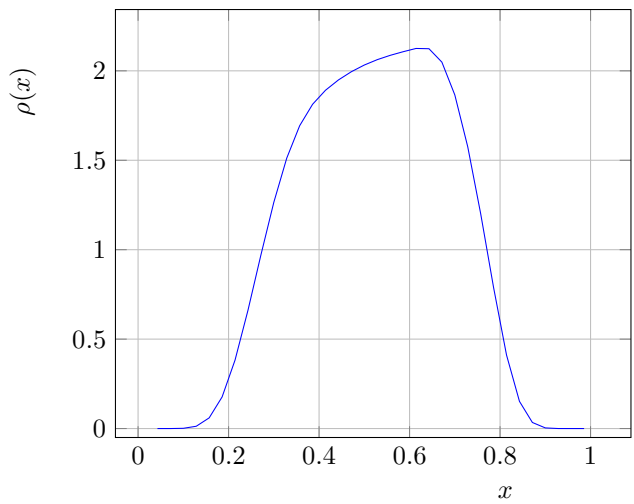
\begin{figure}[h]
\begin{tikzpicture}
\pgfmathsetmacro{\xmin}{-0.05}
\pgfmathsetmacro{\ymin}{-0.05}

\begin{axis}[%ymin=-0.5,
  xmin=\xmin,
  ymin=\ymin,
  ylabel=\textbf{\hskip 100pt $\rho(x)$},
  xlabel=\textbf{\hskip 100pt  $x$},
  grid=major,]
%%%%%%%%%%%%%%%%%%%%%%%%%%%%%%%%%%
\addplot[ color = blue]  table {
   0.042857142857142   0.00000021882838
   0.071428571428571   0.00004875520282
   0.100000000000000   0.00139058027394
   0.128571428571429   0.01285562387977
   0.157142857142857   0.05970418622246
   0.185714285714286   0.17675440458957
   0.214285714285714   0.38266786418593
   0.242857142857143   0.66104987380478
   0.271428571428571   0.96999616545632
   0.300000000000000   1.26494806652392
   0.328571428571429   1.51198047233515
   0.357142857142857   1.69320001992139
   0.385714285714286   1.81307298781228
   0.414285714285714   1.89159301785632
   0.442857142857143   1.94914025550195
   0.471428571428571   1.99517139414845
   0.500000000000000   2.03245061689562
   0.528571428571428   2.06245771429335
   0.557142857142857   2.08659341329150
   0.585714285714286   2.10677687623997
   0.614285714285714   2.12513201303858
   0.642857142857143   2.12365312093870
   0.671428571428571   2.04861852949439
   0.700000000000000   1.86530554055831
   0.728571428571429   1.57532806608034
   0.757142857142857   1.20129796788375
   0.785714285714286   0.78792223052515
   0.814285714285714   0.41103691552378
   0.842857142857143   0.15214214312743
   0.871428571428571   0.03398682564466
   0.900000000000000   0.00360389953375
   0.928571428571429   0.00011973304363
   0.957142857142857   0.00000050732314
   0.985714285714286   0.00000000002044
};
%%%%%%%%%%%%%%%%%%%%%%%%%%%%%%%%%%%%%%%%%%%%%%%%%%%%%%%%%%%%
\end{axis}
\end{tikzpicture}
\caption{The probability density function $\rho(x)$ of the asymptotic distribution of the normalized length of a path towards synchronization. For $N$ sufficiently large and $\delta>0$ sufficiently small, the proportion of paths of length $N(N-1)(x\pm \delta)/2$ is approximatively $\rho(x)\,\delta$.}
\label{fig:DistributionLength}
\end{figure}
\end{center}

\ms Summarizing, the transition diagram for $K_{N}$ is composed by levels ${\rm L}_0,{\rm L}_1,\ldots,{\rm L}_{N(N-1)/2}$, in such a way that each path towards synchronization passes through levels of increasing index until reaching level $N(N-1)/2$ which contains only the complete graph, representing the full $\epsilon$-synchronization. A typical initial condition starting at ${\rm L}_n$, will take $N(N-1)/2-n$ steps to attain the complete graph. The number of subnetworks at level $n=N(N-1)/2-\ell$ is given by $F_N(\ell)$, defined by Equation~\eqref{eq:DistLengths}. The number of subnetworks at each level increases monotonously from $1$ to 
\begin{equation}\label{eq:ModeKn}
{\rm mode}_N(\ell):=\max_{1 \leq \ell \leq N(N-1)/2} F_N(\ell) 
\approx 0.632\,\frac{N\,(N-1)}{2},
\end{equation}
and then decreases monotonously to $1$ as depicted in Figure~\ref{fig:DistributionLength}. Being the distribution of those lengths negatively skewed, the mean length of these paths is smaller than the most frequent length and we have
\[\langle \ell\rangle_N :=\frac{\sum_{\ell=1}^{N(N-1)/2} \ell\,F_N(\ell)}{C_N}
\approx 0.523\,\frac{N\,(N-1)}{2} < {\rm mode}_N(\ell).
\] 

\ms From the calculations above, we can get an idea of some features of a typical synchronization path in the Laplacian of the complete graph, for example, if we were to take a random ordered initial condition of dimension $N$, then its associated synchronization path would most likely be of length as in Equation~(\ref{eq:ModeKn}).

\section{Some results concerning $K_{N,N}$}\label{sec:KNN} 
\ms Let us recall that the Laplacian matrix of $L$ corresponding to the network $K_{N,N}$ has the following entries
\begin{eqnarray*}
L(i,j)   
  &=&\left\{ \begin{array}{rl} 
1,
    & \text{ if } 
    N < i \leq 2 N  \text{ and }  0 < j \leq N 
    \text{ or }
    N < j \leq 2 N  \text{ and }  0 < i \leq N,\\ 
-N,
    & \text{ if } 
    i = j,\,  1\leq i, j \leq 2N,\\
0,
	& \text{ otherwise. } 
    \end{array} \right.
\end{eqnarray*}
An eigenbasis can be computed in terms of the canonical basis and written as the set $\mathcal{B}=\{u^m,v^n,w^n:\ 1\leq m\leq 2,\, 1\leq n\leq N-1\}$, where $u^1=\sum_{k=1}^{2N}\e^k, u^2=\sum_{k=1}^N(\e^k-\e^{k+N})$ and for each $n\geq 1$, $v^n=\e^{n+1}-\e^1$ and $w^n = \e^{N+n+1}-\e^{N+1}$. The Laplacian matrix $L$ acts on this basis as follows: $Lu^1=0$, $Lu^2=-2N\,u^2$ and $Lv^n=-N\,v^n$, $Lw^n=-N\,w^N$ for each $n=1,2,\ldots, N-1$. An initial condition $x\in \IR^{2N}$ can be decomposed as 
\[
x=\bar{x}\,u^1 + (\bar{x}_1-\bar{x})\,u^2+\sum_{n=1}^{N-1}\left((x_{n+1}-\bar{x}_1)\,v^n+(x_{N+n+1}-\bar{x}_2)\,w^n\right),\] 
where 
\begin{equation}\label{eq:MeanCoordinates}
\bar{x}    := \frac{\sum_{n=1}^{2N}x_n}{2N},\,
\bar{x}_1  := \frac{\sum_{n=1}^Nx_n}{N}\, \text{ and }\, 
\bar{x}_2  := \frac{\sum_{n=1}^Nx_{N+n}}{N}.
\end{equation}
Therefore, for all $t\in\IR$ we have
\begin{align*}
x(t)&=\bar{x}\,u^1+ e^{-2N\,t}(\bar{x}_1-\bar{x})\,u^2+e^{-Nt}\sum_{n=1}^{N-1}\left((x_{n+1}-\bar{x}_1),v^n +(x_{N+n+1}-\bar{x}_2)\,w^n\right),\\
    &=\sum_{n=1}^N\left(\left(1-e^{-Nt}\right)\left(\bar{x}-e^{-Nt}\bar{x}_1\right)+e^{-Nt}x_n\right)\e^n \\
    &\hskip 60pt + \sum_{n=1}^N\left(\left(1-e^{-Nt}\right)\left(\bar{x}-e^{-Nt}\bar{x}_2\right)+e^{-Nt}x_{N+n}\right)\e^{N+n}.
\end{align*}
From here it follows that
\begin{eqnarray}
x_n(t)-x_{N+m}(t)     &=&e^{-N\,t}\left(x_n-x_{N+m}+\left(1-e^{-Nt}\right)\left(\bar{x}_1-\bar{x}_2\right)\right),
\label{eq:KNNMonotonicityA} \\ 
x_n(t)-x_m(t)         &=&e^{-N\,t}\left(x_n-x_{m}\right), 
\label{eq:KNNMonotonicityB} \\ 
x_{N+n}(t)-x_{N+m}(t) &=&e^{-N\,t}\left(x_{N+n}-x_{N+m}\right), 
\nonumber
\end{eqnarray}
for all $t\in\IR$ and each $1\leq m,n\leq N$. Hence, the distance between coordinates in the same party of $K_{N,N}$ decreases monotonously, while the distances between coordinates at different parties oscillates at most once, and then decreases to zero. All the differences decreases monotonously if and only if the initial condition satisfies $\bar{x}_1=\bar{x}_2$. In this case the edges $\{n,m\}$ would be included in the synchronized subnetwork $G_{x(t)}$ for all $t\geq t_{n,m}: = \left(\log|x_n-x_{N+m}|-\log(\epsilon)\right)/N$.

\ms Without lost of generality, we may assume that the initial condition is ordered as $x_1\leq x_2\leq\cdots\leq x_N$, $x_{N+1}\leq x_{N+2}\leq\cdots\leq x_{2N}$. By Equation~\eqref{eq:KNNMonotonicityB}  ensures that $x_1(t)\leq x_2(t)\leq\cdots\leq x_N(t)$ and $x_{N+1}(t)\leq x_{N+2}(t)\leq\cdots\leq x_{2N}(t)$ for all $t\in\IR$. We will further assume, when convenient, that $\bar{x}_1=\bar{x}_2$.

\ms Once again, in order to take advantage of the fact that the Laplacian flow preserves the order of the coordinates at each party, we will define the transition diagram not over the synchronized subnetworks but over combinatorial objects that encode the synchronized subnetworks respecting this order. This will simplify the description the transition diagram, mainly in the monotonous case which is achieved when $\bar{x}_1=\bar{x}_2$. We codify the $\epsilon$-synchronized subnetwork $G_x$ defined by $x_1\leq x_2\leq\cdots\leq x_N,\, x_{N+1}\leq x_{N+2}\leq\cdots\leq x_{2N}$, by the couple of functions $\alpha_x,\omega_x:\{1,2,\ldots,N\}\to\{0,1,2,\ldots,N+1\}$ given by
\begin{eqnarray}\label{eq:KNNIncreasingFunctionA}
\alpha_x(n)    
  &=&\left\{ \begin{array}{ll} 
\min\{\ell \leq N:\, x_n-\epsilon\leq x_{N+\ell} \} 
    & \text{ if } 
    x_{2N} \geq x_n-\epsilon,\\ 
N+1 & \text{ if } \
    x_{2N} < x_n-\epsilon, \end{array} \right.\\ \label{eq:KNNIncreasingFunctionB}
\omega_x(n)    
 &=&\left\{ \begin{array}{ll} 
\max\{\ell \leq N:\, x_n+\epsilon \geq x_{N+\ell} \} 
    & \text{ if } \
    x_{N+1} \leq x_n+\epsilon,\\     
0   & \text{ if } \
    x_{N+1} > x_n + \epsilon.\end{array} \right.
\end{eqnarray}
Notice that ${\rm im}(\alpha_x)\subset[1,N+1]$ while ${\rm im}(\omega_x)\subset[0,N]$. Both functions are increasing and such that $\alpha_x(n)\leq \omega_x(n)+1$ for each $1\leq n\leq N$. We present an example of the construction of the increasing functions from a given initial condition, in Figure~\ref{fig:KNN4ConstructionEx}.

\begin{center}
\begin{figure}[h]

\begin{tikzpicture}
%%%%%%%%%% Vecindades
\draw  (-1,7)    -- (1,8.2);
\draw   (-1,7)    -- (1,5.8);
\draw  (-1,3.5)  -- (1,4.7);
\draw  (-1,3.5)  -- (1,2.3);

%%%%%%%%%% Coordinates
\filldraw (-1,7)    circle (3pt);
\filldraw (1,6)     circle (3pt);
\filldraw (1,4)     circle (3pt);
\filldraw (-1,3.5)  circle (3pt);

%%%%%%%%%% Coordinates labels
\draw (1.5,6)    node {\text{$x_4$}};
\draw (1.5,4)    node {\text{$x_3$}};
\draw (-1.5,7)   node {\text{$x_2$}};
\draw (-1.5,3.5) node {\text{$x_1$}};

%%%%%%%%%% Arrows
\draw  (4,6.5) -- (6,6.5);
\draw  (4,3.5) -- (6,3.5);

%%%%%%%%%% Vertices
\filldraw [color=blue] (4,6.5)  circle (3pt);
\filldraw [color=blue] (6,6.5)  circle (3pt);
\filldraw [color=blue] (4,3.5)  circle (3pt);
\filldraw [color=blue] (6,3.5)  circle (3pt);

%%%%%%%%%% Vertex labels
\draw (6.5,6.5) node {\text{4}};
\draw (6.5,3.5) node {\text{3}};
\draw (3.5,6.5) node {\text{2}};
\draw (3.5,3.5) node {\text{1}};

%%%%%%%%%% Text Node
\draw (10,5.5) node {\text{$\alpha_x$=(1,2)}};
\draw (10,4.5) node {\text{$\omega_x$=(1,2)}};

%%%%%%%%%% Label Text 
\draw (0,1)  node {\text{(a)}};
\draw (5,1)  node {\text{(b)}};
\draw (10,1) node {\text{(c)}};
\end{tikzpicture}
\caption{In (a), an example of the relative position of the coordinates of $x=(x_1,x_2,x_3,x_4)$ is illustrated with black dots. The angles that opens from the first two coordinates indicate their $\epsilon$-neighborhood. To construct $G_{x}$, according to Equation~(\ref{eq:NetworkFuctionsCorrespondance}), it is enough to observe that $x_3$ is inside the $\epsilon$-neighborhood of $x_1$, and also $x_4$ is inside the $\epsilon$-neighborhood of $x_2$, hence in (b), vertices 1 and 3 are connected as well as vertices 2 and 4. In (c), the increasing functions determined by $x$ are shown. The function $\alpha_x$ codified the fact that  $x_{3}$ is the first coordinate of the second party inside the angle opening from $x_2$ and
similarly $x_4$ with respect to $x_2$. On the other hand, $\omega_x$, indicates that $x_{3}$ is the last coordinate of the second party inside the angle opening from $x_1$ and respectively $x_4$ with respect to $x_2$.}~\label{fig:KNN4ConstructionEx}
\end{figure}
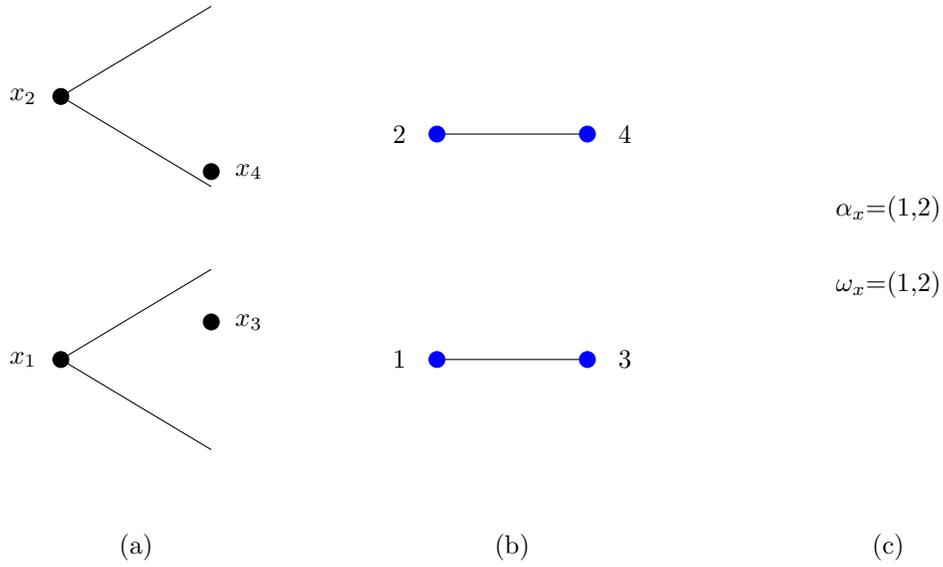
\end{center}

\ms Let $I_N:=\{\phi:\{1,\ldots,N\}\to\{0,\ldots,N+1\}:\, \phi(n+1)\geq \phi(n)\, \text{ for all }\, 1\leq n < N\}$. From the discussion in Appendix~\ref{app:KNNIncreasing}, it follows that the collection
\begin{equation}\label{eq:FunctionsCollection}
\Phi_{N,N}:=\{(\alpha,\omega):\ \alpha,\omega\in I_N: \, {\rm im}(\alpha)\subset[1,N+1],\, {\rm im}(\omega)\subset[0,N] \text{ and }\alpha\leq \omega+1\},
\end{equation}
codify all the $\epsilon$-synchronized subnetworks of $K_{N,N}$ compatible with an ordered initial conditions $x_1\leq x_2\leq\cdots\leq x_N,\, x_{N+1}\leq x_{N+2}\leq\cdots\leq x_{2N}$. The correspondence is given as follows. To $(\alpha,\omega)\in \Phi_{N,N}$ we associate the subnetwork $G_{(\alpha,\omega)}\subset K_{N,N}$ with edges in the set
\begin{equation}~\label{eq:NetworkFuctionsCorrespondance}
E_{(\alpha,\omega)}=\{\{n,N+m\}:\,1\leq n, m\leq N,\text{ and } \alpha(n)\leq m \leq \omega(n)\},
\end{equation}
which is consistent with the fact that $(\alpha,\omega)=(\alpha_x,\omega_x)$ if and only if $G_{(\alpha,\omega)}=G_x$. The correspondence in Equation~\eqref{eq:NetworkFuctionsCorrespondance} establishes a mapping from $\Phi_{N,N}$ to the collection of $\epsilon$-synchronized subnetworks defined by ordered initial conditions, in other words, it is in this case the $\lambda$ mapping associated with Equation~(\ref{eq:transitiondiagram}). The elements in $\Phi_{N,N}$ can be related to combinatorial objects, the parallelo-polyminoes inscribed in a given rectangle. The number of these objects is given by the so called the Narayana numbers~\cite{Stanley1999}. A parallelo-polyminoe in the rectangular lattice of size $p\times q$ is a connected union of squares delimited by two increasing boundary functions $L,U:\{1,2,\ldots,p\}\to\{0,1,\ldots,q\}$ such that $L(1)=0$, $U(p)=q$, and $L(n) < U(n-1)$ for each $2\leq n\leq p$.
\begin{center}
\begin{figure}[h]

\begin{tikzpicture}[scale=0.7]
\draw[help lines] (1,0) grid (15,7);
%%%%%%%%%%%%%%%%%%%%%%%%%%%%%%%%%%%%%%%
%%%%%%%%%%%%%%%%%%%%%%%%%%%%%%%%%%%%%%%%%%
\draw[blue, thick] (1,0)--(6,0);
\draw[blue, thick] (6,0)--(6,2);
\draw[blue, thick] (6,2)--(10,2);
\draw[blue, thick] (10,2)--(10,5);
\draw[blue, thick] (10,5)--(13,5);
\draw[blue, thick] (13,5)--(15,5);
\draw[blue, thick] (15,5)--(15,7);
%%%%%%%%%%%%%%%%%%%%%%%%%%%%%%%%%%%%%%%%%
\draw[red, thick] (1,0)--(1,1);
\draw[red, thick] (1,1)--(4,1);
\draw[red, thick] (4,1)--(4,3);
\draw[red, thick] (4,3)--(7,3);
\draw[red, thick] (7,3)--(7,5);
\draw[red, thick] (7,5)--(9,5);
\draw[red, thick] (9,5)--(9,6);
\draw[red, thick] (9,6)--(13,6);
\draw[red, thick] (13,6)--(13,7);
\draw[red, thick] (13,7)--(15,7);
%%%%%%%%%%%%%%%%%%%%%%%%%%%%%%%%%%%%%%%%%
\end{tikzpicture}
\caption{A parallelo-polyminoe in the lattice of size $14\times 10$. The blue path defines the lower border function $L=(0,0,0,0,0,2,2,2,2,5,5,5,5,5)$, while the red one defines the upper border $U=(1,1,1,3,3,3,5,5,6,6,6,6,7,7)$.}~\label{fig:Parallelo-polyminoe}
\end{figure}
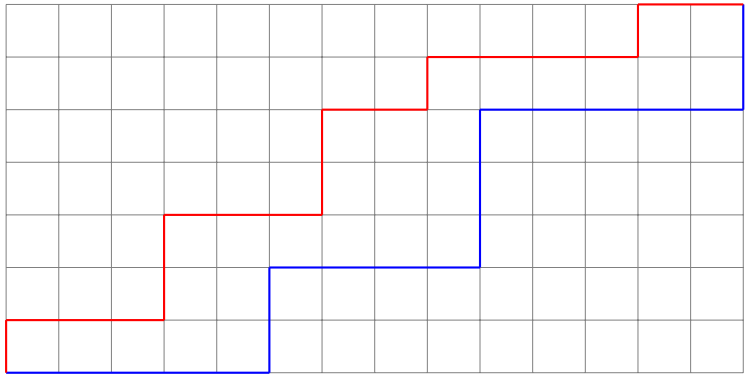 
\end{center}
The number of parallelo-polyminoes in the lattice of size $p\times q$ 
is given by the Narayana number~\cite{BarcucciFrosiniRinaldi2005}
\begin{equation}\label{eq:NarayanaNumbers}
T(p+q-1,q):=\frac{1}{p+q-1}
\left(\begin{matrix}p+q-1\\q\end{matrix}\right)
\left(\begin{matrix}p+q-1\\q-1\end{matrix}\right).
\end{equation}
To each couple $(\alpha,\omega)\in\Phi_{N,N}$ we associate a parallelo-polyminoe in $\{0,1,\ldots,N+1\}\times \{0,1,\ldots, N+1\}$ with border functions $L,U:\{1,\ldots,N+1\}\to\{0,1,\ldots,N+1\}$, such that
\begin{equation}~\label{eq:Polyminoe}
L(n)=\left\{\begin{array}{cl} 
            0  & \text{ for } n = 1,\\
      \alpha(n-1)-1 & \text{ for } 2\leq n \leq N+1,\\
            \end{array}\right.
\text{ and } \
U(n)=\left\{\begin{array}{cl}
     \omega(n)+1 & \text{ for } 1\leq n \leq N,\\
            N+1  & \text{ for } n=N+1.\end{array}\right.
\end{equation}
In this way, we establish a one-to-one correspondence between parallelo-polyminoes and couples in $\Phi_{N,N}$, from which we obtain
\begin{equation}\label{eq:CardPhiNN}
|\Phi_{N,N}|=T(2N+1,N+1)=\frac{1}{2N+1}
\left(\begin{matrix}2N+1\\N+1\end{matrix}\right)
\left(\begin{matrix}2N+1\\N\end{matrix}\right).
\end{equation} 

\ms Thanks to the equivalence given by the Equation~\eqref{eq:KNNetworkFuctionCorrespondance}, each sequence of $\epsilon$-synchronized subnetworks defined by an ordered initial condition is faithfully codified by the corresponding sequences of couples of increasing functions given by the Equations~\eqref{eq:KNNIncreasingFunctionA} and~\eqref{eq:KNNIncreasingFunctionB}. As mentioned above, for an initial condition $x\in\IR^{2N}$ such that $\bar{x}_1=\bar{x}_2$, all the differences $x_{N+m}(t)-x_n(t)$ converge to $0$ monotonously and at the same speed. We will say that such initial conditions are balanced. In this case, each one of the maps $t\mapsto \alpha_{x(t)}$ and $t\mapsto \omega_{x(t)}$ are coordinate-wise monotonous, and they converge respectively to the constant functions ${\bf 1}(n)=1$ and ${\bf N}(n)=N$ at time $t_{1,N}=\left(\log|x_1-x_{2N}|-\log(\epsilon)\right)/N$. The sequence of switching times $0 < t_1 < t_2 < \cdots < t_\ell$ is such that $(\alpha_{x(t_\tau)},\omega_{x(t_\tau)})\neq (\alpha_{x(t_{\tau+1})},\omega_{x(t_{\tau+1})})$. Let us denote $\alpha_{t_\tau}$ by $\alpha_\tau$, and the corresponding for $\omega$. For a typical initial condition, at each switching time only one of the functions $\alpha_\tau$ or $\omega_\tau$ changes and it changes only at one site. The sequence $((\alpha_0,\omega_0),(\alpha_1,\omega_1),\ldots,(\alpha_\ell,\omega_\ell))$ can be  determined by the initial couple $(\alpha_0,\omega_0)$, the jump sites $n_1,n_2,\ldots,n_\ell\in\{1,2,\ldots,N\}^\ell$ and binary labels $q_1,q_2,\ldots,q_\ell\in(-1,+1)^\ell$ as follows:
\begin{equation}\label{eq:EvolutionAlphaOmega}
(\alpha_{\tau+1},\omega_{\tau+1})=\left\{
\begin{array}{ll}
(\alpha_\tau-\delta_{n_\tau}, \omega_\tau) 
                          & \text{ if } q_\tau=-1,\\
(\alpha_\tau, \omega_\tau+\delta_{n_\tau}) 
                          & \text{ if } q_\tau=+1.
\end{array}\right.
\end{equation}
To the couple $(\alpha_\tau,\omega_\tau)$, we can associate a parallelo-polyminoe according to  Equation~\eqref{eq:Polyminoe}. In the transition $(\alpha_\tau,\omega_\tau)\to(\alpha_{\tau+1},\omega_{\tau+1})$, the area inside the corresponding parallelo-polyminoe increases by one unit until the final area $N\times N$. 

\ms Realizable sequences $((n_1,q_1),(n_2,q_2),\ldots,(n_\ell,q_\ell))$, are those compatible with a balanced initial condition $x\in\IR^{2N}$ and are completely determined by the differences $\Delta_{n,m}:=x_{N+m}-x_n$ with $1\leq n,m\leq N$ as follows: For $\epsilon < |\Delta_{n_1,m_1}| < |\Delta_{n_2,m_2}| < \cdots < |\Delta_{n_{N^2},m_{N^ 2}}|$ we have the sequence $((n_1,q_1),(n_2,q_2),\ldots,(n_{N^2},q_{N^2}))$, where $q_\tau={\rm sign}(\Delta_{n_\tau,m_\tau)}$ for each $1\leq \tau \leq N^2$. If we consider all the possible orderings $\Delta:=\{\Delta_{n,m}:\, 1\leq n,m\leq N\}$
compatible with an initial condition, not necessarily balanced, and we assume that the dynamics towards synchronization is completely determined by this ordering as in the balanced case, we obtain a transition diagram with vertices in $\Phi_{N,N}$ with maximal paths starting at the couples $(\alpha,\omega)$ codifying the disconnected subnetwork, and ending at the couple $({\bf 1},{\bf N})$ which codifies $K_{N,N}$. This digraph contains all the paths towards synchronization starting at balanced initial conditions but it also contains paths which are not compatible with any balanced initial condition. For instance, in the case $N=2$ there are 20 realizable possible orderings $\{\Delta_{n,m}:\, 1\leq n,m\leq N\}$, which we depict in Table~\ref{tab:NN2Ordering}, defining 20 paths towards synchronization represented in the transition diagram of Figure~\ref{fig:N2TransitionGraph}. Nevertheless, there are 4 orderings, and therefore 4 paths towards synchronization, which are incompatible with a balanced initial condition. The coordinate arrangements incompatible with a balanced initial conditions are $x_1<x_2<x_3<x_4$ and $x_3<x_4<x_1<x_2$. In general there are 2 arrangements of initial conditions, $x_1<\cdots <x_N<x_{N+1}<\cdots<x_{2N}$ and $x_{N+1}<\cdots <x_{2N}<x_1<\cdots<x_N$, which are incompatible with a balanced initial condition. These arrangements define maximal paths starting at vertices $({\bf 1},{\bf 0})$ and $({\bf N+1},{\bf N})$, which for the case $N=2$ we indicate in red in Figure~\ref{fig:N2TransitionGraph}. 

\begin{table}[h]
\begin{tabular}{|c|c|c|}
\hline 
Coordinates    &   Differences       & Signs \\
\hline
$x_1<x_2<x_3<x_4$ & 
      $|\Delta_{2,1}|<|\Delta_{2,2}|<|\Delta_{1,1}|<|\Delta_{1,2}|$ 
                                         & $(+1,+1,+1,+1)$ \\ 
                  &
      $|\Delta_{2,1}|<|\Delta_{1,1}|<|\Delta_{2,2}|<|\Delta_{1,2}|$ 
                                         & $(+1,+1,+1,+1)$ \\ 
\hline
$x_1<x_3<x_2<x_4$ & 
      $|\Delta_{2,1}|<|\Delta_{2,2}|<|\Delta_{1,1}|<|\Delta_{1,2}|$ 
                                         & $(-1,+1,+1,+1)$ \\ 
                  & 
      $|\Delta_{2,2}|<|\Delta_{2,1}|<|\Delta_{1,1}|<|\Delta_{1,2}|$ 
                                         & $(+1,-1,+1,+1)$ \\ 
                  & 
      $|\Delta_{2,2}|<|\Delta_{1,1}|<|\Delta_{2,1}|<|\Delta_{1,2}|$ 
                                         & $(+1,+1,-1,+1)$ \\ 
                  & 
      $|\Delta_{2,1}|<|\Delta_{1,1}|<|\Delta_{2,2}|<|\Delta_{1,2}|$ 
                                         & $(-1,+1,+1,+1)$ \\ 
                  & 
      $|\Delta_{1,1}|<|\Delta_{2,1}|<|\Delta_{2,2}|<|\Delta_{1,2}|$ 
                                         & $(+1,-1,+1,+1)$ \\ 
                  & 
      $|\Delta_{1,1}|<|\Delta_{2,2}|<|\Delta_{2,1}|<|\Delta_{1,2}|$ 
                                         & $(+1,+1,-1,+1)$ \\ 
\hline
$x_1<x_3<x_4<x_2$ & 
      $|\Delta_{1,1}|<|\Delta_{2,2}|<|\Delta_{1,2}|<|\Delta_{2,1}|$ 
                                         & $(+1,-1,+1,-1)$\\ 
                  & 
      $|\Delta_{2,2}|<|\Delta_{1,1}|<|\Delta_{2,1}|<|\Delta_{1,2}|$ 
                                         & $(-1,+1,-1,+1)$\\ 
\hline
$x_3<x_4<x_1<x_2$ & 
      $|\Delta_{1,2}|<|\Delta_{1,1}|<|\Delta_{2,2}|<|\Delta_{2,1}|$ 
                                         & $(-1,-1,-1,-1)$\\ 
                  & 
      $|\Delta_{1,2}|<|\Delta_{2,2}|<|\Delta_{1,1}|<|\Delta_{2,1}|$ 
                                         & $(-1,-1,-1,-1)$\\ 
\hline
$x_3<x_1<x_4<x_2$ & 
      $|\Delta_{1,2}|<|\Delta_{1,1}|<|\Delta_{2,2}|<|\Delta_{2,1}|$ 
                                         & $(+1,-1,-1,-1)$ \\ 
                  & 
      $|\Delta_{1,1}|<|\Delta_{1,2}|<|\Delta_{2,2}|<|\Delta_{2,1}|$ 
                                         & $(-1,+1,-1,-1)$ \\ 
                  & 
      $|\Delta_{1,1}|<|\Delta_{2,2}|<|\Delta_{1,2}|<|\Delta_{2,1}|$ 
                                         & $(-1,-1,+1,-1)$ \\ 
                  & 
      $|\Delta_{1,2}|<|\Delta_{2,2}|<|\Delta_{1,1}|<|\Delta_{2,1}|$ 
                                         & $(+1,-1,-1,-1)$ \\ 
                  & 
      $|\Delta_{2,2}|<|\Delta_{1,2}|<|\Delta_{1,1}|<|\Delta_{2,1}|$ 
                                         & $(-1,+1,-1,-1)$ \\
                  & 
     $|\Delta_{2,2}|<|\Delta_{1,1}|<|\Delta_{1,2}|<|\Delta_{2,1}|$ 
                                         & $(-1,-1,+1,-1)$ \\ 
\hline
$x_3<x_1<x_2<x_4$ & 
     $|\Delta_{1,1}|<|\Delta_{2,2}|<|\Delta_{1,2}|<|\Delta_{2,1}|$ 
                                         & $(-1,+1,+1,+1)$ \\ 
                  & 
     $|\Delta_{2,2}|<|\Delta_{1,1}|<|\Delta_{2,1}|<|\Delta_{1,2}|$ 
                                         & $(+1,-1,-1,+1)$ \\
\hline
\end{tabular}

\ms
\caption{The twenty different orderings of the differences between coordinates at opposite parties, and corresponding signs, for a typical initial conditions in $\IR^4$.}~\label{tab:NN2Ordering}
\end{table}

\begin{center}
\begin{figure}[h!]

\begin{tikzpicture}[scale=2]

%%%%%%%%%% Arrows
\draw [red, thick,-to]  (1,8) -- (1,7.1);
\draw [-to] 			(1,7) -- (1,6.1);
\draw [-to] 			(1,6) -- (1,5.1); 
\draw [-to]				(2,8) -- (2,7.1);
\draw [-to] 			(2,7) -- (2,6.1);
\draw [-to] 			(2,8) -- (1.1,7.1);
\draw [-to] 			(2,7) -- (1.1,6.1);
\draw [-to] 			(2,6) -- (1.1,5.1);
\draw [-to] 			(2,5) -- (1.1,5);
\draw [-to] 			(1,5) -- (1,4) -- (2.9,4);
\draw [-to] 			(4,1) -- (4,1.9);
\draw [-to] 			(3,2) -- (3,2.9);
\draw [-to] 			(4,1) -- (3.1,1.9);
\draw [-to] 			(4,2) -- (3.1,2.9);
\draw [red, thick, -to] (3,1) -- (3,1.9);
\draw [-to] 			(3,2) -- (3.9,2.9);
\draw [-to] 			(4,3) -- (4,3.9); 
\draw [-to] 			(3,3) -- (3.9,3.9);
\draw [-to] 			(4,8) -- (4.9,7.1);
\draw [-to] 			(5,7) -- (5,5.1);
\draw [-to] 			(5,5) -- (4.1,4.1); 
\draw [-to] 			(4,4) -- (3.1,4); 
\draw [-to] 			(2,8) -- (2.9,7.1);
\draw [-to] 			(3,7) -- (2.1,6.1);
\draw [-to] 			(3,7) -- (3,6) -- (2.1,5.1);
\draw [-to] 			(6,6) -- (5.1,5.1);
\draw [-to] 			(6,4) -- (5.1,4.9);
\draw [-to] 			(6,5) -- (6,5.9);
\draw [-to] 			(6,5) -- (6,4.1);
\draw [-to] 			(4,8) -- (2.1,7.1);
\draw [-to] 			(4,2) -- (5,3) -- (5,4.9);
\draw [-to] 			(1,7) -- (1.9,5.1);
%%%%%%%%%%%%%%%%%%%%%%%%%%%%%%%%%%%%%%%%%%%%%%%

%%%%%%%%%%%%% Vertices and their labels
\filldraw [color=blue] (1,8) circle (1.25pt);
\draw (1,8.25) node {\rojo\text{\underline{(11,00)}}};
%%%%%%%%%%%%%%%%%%%%%%%%%%%%%%%%%%%%%%%%%%%%%%%%
\filldraw [color=blue] (1,7) circle (1.25pt);
\draw (0.6,7) node {\text{(11,01)}};
%%%%%%%%%%%%%%%%%%%%%%%%%%%%%%%%%%%%%%%%%%%%%%%%
\filldraw [color=blue] (1,6) circle (1.25pt);
\draw (0.6,6) node {\text{(11,02)}};
%%%%%%%%%%%%%%%%%%%%%%%%%%%%%%%%%%%%%%%%%%%%%%%
\filldraw [color=blue] (1,5) circle (1.25pt);
\draw (0.6,5) node {\text{(11,12)}};
%%%%%%%%%%%%%%%%%%%%%%%%%%%%%%%%%%%%%%%%%%%%%%%
\draw (2,8.25) node {\text{\underline{(12,01)}}};
\filldraw [color=blue] (3,1) circle (1.25pt);
%%%%%%%%%%%%%%%%%%%%%%%%%%%%%%%%%%%%%%%%%%%%%%%
\filldraw [color=blue] (3,2) circle (1.25pt);
\draw (2.6,2) node {\text{(23,22)}};
%%%%%%%%%%%%%%%%%%%%%%%%%%%%%%%%%%%%%%%%%%%%%%%
\filldraw [color=blue] (3,3) circle (1.25pt);
\draw (2.6,3) node {\text{(13,22)}};
%%%%%%%%%%%%%%%%%%%%%%%%%%%%%%%%%%%%%%%%%%%%%%
\filldraw [color=blue] (2,5) circle (1.25pt);
\draw (2.4,5) node {\text{(11,11)}};
%%%%%%%%%%%%%%%%%%%%%%%%%%%%%%%%%%%%%%%%%%%%%%
\filldraw [color=blue] (2,6) circle (1.25pt);
\draw (2.4,6) node {\text{(12,12)}};
%\node[rectangle,draw,fill=white] (r) at (285,140) {3};
%%%%%%%%%%%%%%%%%%%%%%%%%%%%%%%%%%%%%%%%%%%%%%%
\filldraw [color=blue] (2,7) circle (1.25pt);
\draw (1.6,7) node {\text{(12,02)}};
%%%%%%%%%%%%%%%%%%%%%%%%%%%%%%%%%%%%%%%%%%%%%
\filldraw [color=blue] (2,8) circle (1.25pt);
%%%%%%%%%%%%%%%%%%%%%%%%%%%%%%%%%%%%%%%%%%%%%%%
\draw (3,0.75) node {\rojo\text{\underline{(33,22)}}};
%%%%%%%%%%%%%%%%%%%%%%%%%%%%%%%%%%%%%%%%%%%%%%%
\filldraw [color=blue] (4,1) circle (1.25pt);
\draw (4,0.75) node {\text{\underline{(23,12)}}};
%%%%%%%%%%%%%%%%%%%%%%%%%%%%%%%%%%%%%%%%%%%%%%%
\filldraw [color=blue] (4,2) circle (1.25pt);
\draw (4.4,2) node {\text{(13,12)}};
%%%%%%%%%%%%%%%%%%%%%%%%%%%%%%%%%%%%%%%%%%%%%%%
\filldraw [color=blue] (4,3) circle (1.25pt);
\draw (4.4,3) node {\text{(22,22)}};
%%%%%%%%%%%%%%%%%%%%%%%%%%%%%%%%%%%%%%%%%%%%%%%
\filldraw [color=blue] (3,4) circle (1.25pt);
\node[rectangle, draw] (r) at (3,4.3) {\text{(11,22)}};
%%%%%%%%%%%%%%%%%%%%%%%%%%%%%%%%%%%%%%%%%%%%%%%%
\filldraw [color=blue] (3,7) circle (1.25pt);
\draw (3.4,7) node {\text{(12,11)}};
%%%%%%%%%%%%%%%%%%%%%%%%%%%%%%%%%%%%%%%%%%%%%%%%
\draw (4,8.25) node {\text{\underline{(13,02)}}};
\filldraw [color=blue] (4,8) circle (1.25pt);
%%%%%%%%%%%%%%%%%%%%%%%%%%%%%%%%%%%%%%%%%%%%%%%%
\filldraw [color=blue] (5,7) circle (1.25pt);
\draw (5.4,7) node {\text{(13,12)}};
%%%%%%%%%%%%%%%%%%%%%%%%%%%%%%%%%%%%%%%%%%%%%%%%%
\filldraw [color=blue] (5,5) circle (1.25pt);
\draw (4.6,5) node {\text{(12,12)}};
%%%%%%%%%%%%%%%%%%%%%%%%%%%%%%%%%%%%%%%%%%%%%%%%%%
\filldraw [color=blue] (4,4) circle (1.25pt);
\draw (4.4,4) node {\text{(22,22)}};
%%%%%%%%%%%%%%%%%%%%%%%%%%%%%%%%%%%%%%%%%%%%%%%%%%%
\filldraw [color=blue] (6,6) circle (1.25pt);
\draw (6,6.2) node {\text{(12,11)}};
%%%%%%%%%%%%%%%%%%%%%%%%%%%%%%%%%%%%%%%%%%%%%%%%%%%%
\filldraw [color=blue] (6,5) circle (1.25pt);
\draw (6.5,5) node {\text{\underline{(22,11)}}};
%%%%%%%%%%%%%%%%%%%%%%%%%%%%%%%%%%%%%%%%%%%%%%%%%%%%
\filldraw [color=blue] (6,4) circle (1.25pt);
\draw (6,3.8) node {\text{(22,12)}};
%%%%%%%%%%%%%%%%%%%%%%%%%%%%%%%%%%%%%%%%%%%%%%%%%%%%%
\end{tikzpicture}

\caption{The transition diagram which contains all the paths towards synchronization of the Laplacian dynamics on $K_{2,2}$. Each one of the functions $\alpha,\omega$, are codified by a two-digit string. There are six starting configurations, underlined in the diagram, all of them coding the disconnected network. The ending vertex, $(11,22)$, is the couple codifying the complete bipartite graph $K_{2,2}$. In red we indicate the starting couples which are incompatible with a balanced initial condition. In this case, by erasing the elements in color red, we obtain the transition diagram codifying all the paths towards synchronization for balanced initial conditions.}~\label{fig:N2TransitionGraph}
\end{figure}
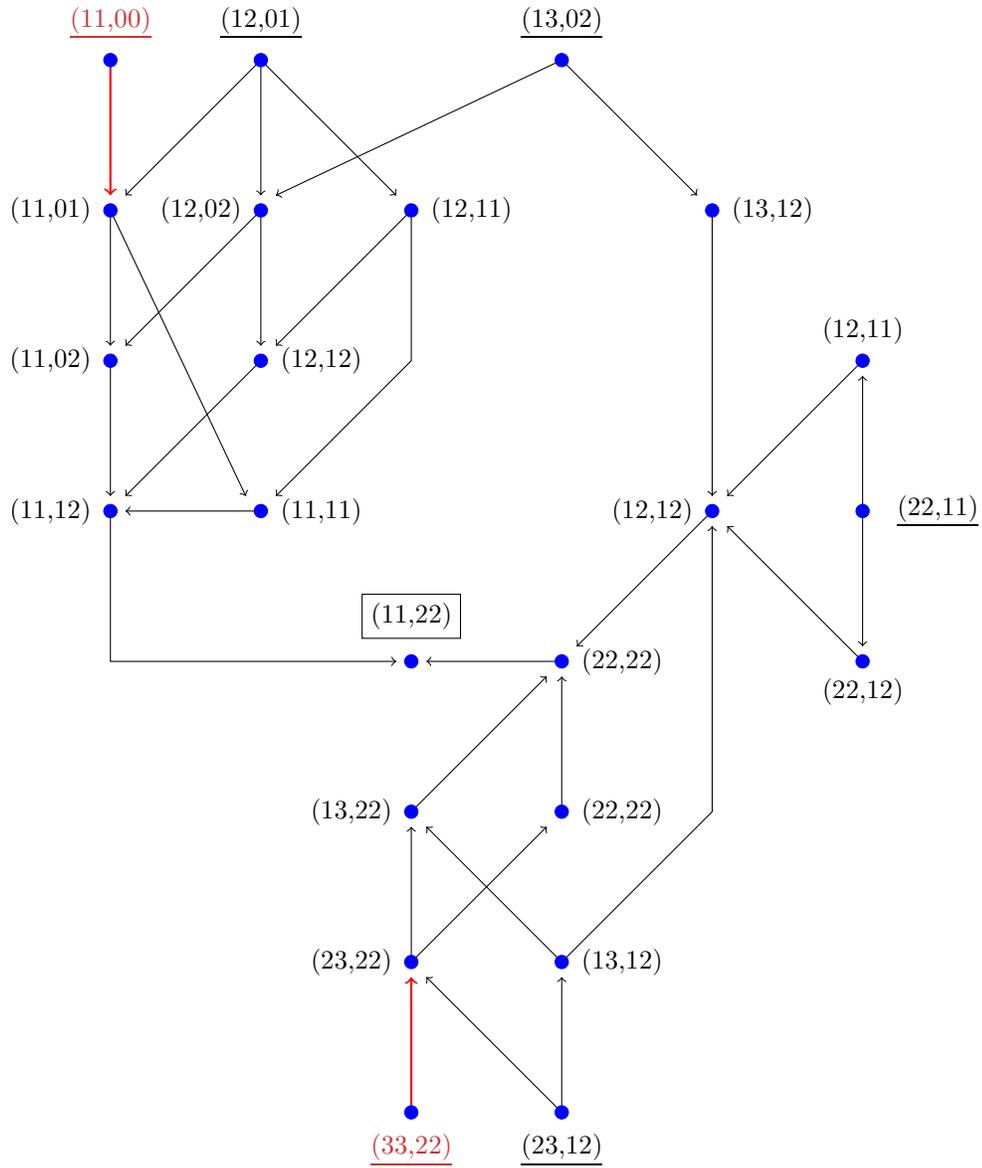 
\end{center}

\ms An easy upper bound for the number of paths towards synchronization starting at typical balanced initial conditions is the following. For each one of the arrangements $x_{i_1}<x_{i_2}<\cdots<x_{i_{2N}}$, obtaining by inter-placing the first $N$ coordinates with respect to the last $N$ coordinates while maintaining the order inside each group of coordinates, there are Golomb$(2N)$ different orderings for the differences $x_{i_k}-x_{i_\ell}$. Each one of these orderings give place to a path towards synchronization, but this path does not depend on the differences between coordinates of the same group (first $N$ or last $N$ coordinates). Furthermore, there are two coordinate arrangements which are incompatible with a balanced initial condition, when $x_{1}<x_{2}<\cdots<x_{2N}$ and when $x_{N+1}<x_{N+2}<\cdots<x_{2N}<x_1<x_{2}<\cdots<x_{N}$, hence the number of paths towards synchronization is upper bounded by
\begin{equation}\label{eq:NoPathsKNN}
\text{Number of paths towards synchronization for }K_{N,N} \leq \left(\left(\begin{matrix} 2N\\N\end{matrix}\right)-2\right)\text{Golomb}(2N).
\end{equation}
As mentioned above, the growth of this quantity with respect to $N$ defines a complexity function analogous to the topological complexity as a function of time.

\ms Similar to the case $K_N$, the number of paths towards synchronization of a given length, $F_{N,N}(\ell)$, is given by the number of couples $(\alpha,\omega)\in\Phi_{N,N}$ such that the corresponding parallelo-polyminoe has an interior with area of $(N+1)^2-\ell$ units. Hence, 
\begin{equation}\label{eq:DistLengthsNN}
F_{N,N}(\ell):=\left|\left\{(\alpha,\omega)\in\Phi_{N,N}:\ \sum_{n=1}^{N+1}(U(n)-L(n))=(N+1)^2-\ell\right\}\right|.
\end{equation}
Here, $L,U:\{1,\ldots,N+1\}\to\{0,1,\ldots,N+1\}$ are the polyminoe border functions defined from the couple $(\alpha,\omega)$ by Equations~\eqref{eq:Polyminoe}. Table~\ref{tab:FNN} shows the distributions $F_{N,N}(\ell)$ for $2\leq N\leq 7$.

\ms
\begin{table}[h!]
\centering
\begin{tabular}{|c|l|}
\hline
$N$ & $F_{N,N}(\ell)$                 \\ \hline
2 & (1,2,5,6,6)                         \\ \hline
3 & (1,2,5,10,16,24,31,36,30,20)        \\ \hline
4 & (1,2,5,10,20,32,53,78,111,146,187,216,243,240,210,140,70) \\ \hline
5 & (1,2,5,10,20,36,61,98,153,228,327,454,611,798,1005,1236,1466,1688,1862,1980,1971,\\  
  & 1850,1540,1120,630,252)\\ \hline
6 & (1,2,5,10,20,36,65,106,173,268,409,600,867,1212,1671,2244,2966,3826,4868,6056,\\   
  & 7422,8906,10519,12166,13830,15352,16704,17656,18133,17890,16903,14966,12306,\\
  & 8988,5670,2772,924)\\ \hline
7 & (1,2,5,10,20,36,65,110,181,288,449,680,1013,1474,2107,2958,4088,5558,7450,9842,\\
  & 12820,16488,20932,26246,32507,39790,48116,57538,67984,79414,91653,104578,117806,\\  
  & 131096,143865,155692,165779,173530,177877,178282,173616,163632,147855,127092,  \\ 
  & 102060,75432,49434,27720,12012,3432)\\\hline
8 & (1,2,5,10,20,36,65,110,185,296,469,720,1093,1618,2369,3400,4824,6732,9296,12654,\\
  & 17054,22694,29912,38976,50333,64320,81489,102242,127219, 156850,191841,232602, \\
  & 279832,333830,395204,464030,540737,625028,716966,815766,920990,1031168,1145253,\\ 
  & 1260882,1376172,1487820,1593022,1687242,1766791,1826112,1860845,1865122, \\
  & 1834995,1765746,1656541,1506540,1320987,1106748,877470,647592,437118,260832, \\
  & 132132,51480,12870)
 \\ \hline
\end{tabular}

\ms
\caption{Number $F_{N,N}(\ell)$ of couples $(\alpha,\omega)\in\Phi_{N,N}$ codifying a subnetworks starting a synchronizing path of length $\ell$.}~\label{tab:FNN}
\end{table}

\ms For each $N$ and $0\leq \ell \leq N$,  the integer $F_{N,N}(\ell)$ coincides with the $\ell$-th term of the Sloans's sequence (Entry A000712 of the On-line Encyclopedia of Integer Sequences~\cite{web:oeisA000712}), which among other things, counts the number of couples of integer partitions $P=(p_1\geq p_2\geq\cdots\geq p_k)$, $Q=(q_1\geq q_2\geq \cdots \geq q_r)$, such that $\sum_{i=1}^kp_i+\sum_{j=1}^rq_j=\ell$. Indeed, we can associate to each such couple of integer partitions $(P, Q)$, a unique couple $L, U:\{1,2,\ldots,N+1\}\to\{0,1,\ldots,N+1\}$ of upper and lower border functions such that $U(i)=N+1-p_i$ and $L(N+2-j)=q_j$. Clearly $\sum_{i=1}^kp_i+\sum_{j=1}^rq_j=\ell$ if and only if the area of the parallelo-polyminoe with border functions $L$ and $U$ is $(N+1)^2-\ell$. The correspondence between integer partitions and border functions cannot go further than $\ell=N$, since for $\ell=N+1$ the couple $((N+1),(0))$ of partitions does not define admissible border functions. On the opposite extreme, $F_{N,N}(N^2)$ counts all the parallelo-polyminoes in $\{0,1,\ldots,N+1\}\times\{0,1,\ldots,N+1\}$ composed of $2N+1$ squares. These squares are arranged in a path going from $(0,0)$ to $(N+1,N+1)$, the next square place at the left or on top of the previous one. Each one of these arrangements can therefore by codified a sequences $(a_1,a_2,\ldots,a_{2N})\in\{L,T\}^{2N}$, with exactly $N$ entries equal to $T$. From this it follows that 
\begin{equation}\label{eq:FNNN2}
F_{N,N}(N^2)=\left(\begin{matrix}2N\\N \end{matrix}\right)
\end{equation} 

\ms The normalized cumulative distribution, $f_{N,N}:[0,1]\to[0,1]$, is given by
\begin{equation}~\label{eq:ComulativeDistNN}
f_{N,N}(x)=\frac{1}{|\Phi_{N,N}|}
\sum_{n\leq x\times N^ 2}F_{N,N}(x),
\end{equation}
where $F_{N,N}$ is given by Equation~\eqref{eq:DistLengthsNN} and $|\Phi_{N,N}|$ by Equation~\eqref{eq:CardPhiNN}. We numerically computed $f_{N,N}(x)$ for increasing $N$, and observe how it approaches a limit distribution $x\mapsto {\rm f}(x)$ whose density $\varrho(x):=d\,{\rm f}(x)/dx$ approaches  the curve depicted in Figure~\ref{fig:DistributionLengthNN}. As for $K_N$, our numerical computation suggest that $\varrho$ is continuous, unimodal, and negatively skewed. 

\begin{center}
\begin{figure}[h!]
\begin{tikzpicture}
\pgfmathsetmacro{\xmin}{-0.05}
\pgfmathsetmacro{\ymin}{-0.05}

\begin{axis}[%ymin=-0.5,
  xmin=\xmin,
  ymin=\ymin,
  ylabel=\textbf{\hskip 100pt $\varrho(x)$},
  xlabel=\textbf{\hskip 100pt  $x$},
  grid=major,]
%%%%%%%%%%%%%%%%%%%%%%%%%%%%%%%%%%
\addplot[ color = blue]  table {
   0.000000000000000   0.000001841021997
   0.015625000000000   0.000003682043995
   0.031250000000000   0.000009205109987
   0.046875000000000   0.000018410219973
   0.062500000000000   0.000036820439947
   0.078125000000000   0.000066276791904
   0.093750000000000   0.000119666429827
   0.109375000000000   0.000202512419707
   0.125000000000000   0.000340589069507
   0.140625000000000   0.000544942511212
   0.156250000000000   0.000863439316751
   0.171875000000000   0.001325535838082
   0.187500000000000   0.002012237043089
   0.203125000000000   0.002978773591690
   0.218750000000000   0.004361381111690
   0.234375000000000   0.006259474790943
   0.250000000000000   0.008881090115150
   0.265625000000000   0.012393760086068
   0.281250000000000   0.017114140487238
   0.296875000000000   0.023296292354293
   0.312500000000000   0.031396789142573
   0.328125000000000   0.041780153207549
   0.343750000000000   0.055068649984323
   0.359375000000000   0.071755673368179
   0.375000000000000   0.092664160191927
   0.390625000000000   0.118414534868669
   0.406250000000000   0.150023041540935
   0.421875000000000   0.188229771051655
   0.437500000000000   0.234212977479123
   0.453125000000000   0.288764300282194
   0.468750000000000   0.353183500990988
   0.484375000000000   0.428225398624411
   0.500000000000000   0.515176867558603
   0.515625000000000   0.614588373370767
   0.531250000000000   0.727579257435284
   0.546875000000000   0.854289437423950
   0.562500000000000   0.995508711773623
   0.578125000000000   1.150690296951095
   0.593750000000000   1.319950177342197
   0.609375000000000   1.501843150679021
   0.625000000000000   1.695562849326733
   0.640625000000000   1.898402970949248
   0.656250000000000   2.108435965515357
   0.671875000000000   2.321311498045353
   0.687500000000000   2.533562924118251
   0.703125000000000   2.739109348076851
   0.718750000000000   2.932788544240621
   0.734375000000000   3.106249636829645
   0.750000000000000   3.252701095695748
   0.765625000000000   3.361912361599733
   0.781250000000000   3.425856578633214
   0.796875000000000   3.433730629715821
   0.812500000000000   3.378266160002071
   0.828125000000000   3.250777227708532
   0.843750000000000   3.049728420489424
   0.859375000000000   2.773573279866986
   0.875000000000000   2.431966125195249
   0.890625000000000   2.037547413507923
   0.906250000000000   1.615441572002658
   0.921875000000000   1.192231117298991
   0.937500000000000   0.804743853431636
   0.953125000000000   0.480197449609214
   0.968750000000000   0.243257918552036
   0.984375000000000   0.094775812422871
   1.000000000000000   0.023693953105718
   };
%%%%%%%%%%%%%%%%%%%%%%%%%%%%%%%%%%%%%%%%%%%%%%%%%%%%%%%%%%%%
%%%%%%%%%%%%%%%%%%%%%%%%%%%%%%%%%%%%%%%%%%%%%%%%%%%%%%%%%%%%
\end{axis}
\end{tikzpicture}
\caption{The probability density function $\varrho(x)$ of the asymptotic distribution of the normalized length of a path towards synchronization. For $N$ sufficiently large and $\delta>0$ sufficiently small, the number of paths of length $N^2(x\pm \delta)/2$ is approximatively $\varrho(x)\,\delta$.}
\label{fig:DistributionLengthNN}
\end{figure}
\end{center}

\ms As we have already mentioned, in the case of $K_{N,N}$ we do not have the complete panorama of its paths towards synchronization, since our methodology is limited to the initial conditions that are balanced. In addition, currently there are no results in combinatorics that allow us to make calculations for arbitrarily large sizes. Nevertheless by directly computing these distributions for low dimensions, we observe a very fast convergence of the normalized distribution $f_{N,N}$. We obtain a unimodal distribution with maximum at  
\begin{equation}\label{eq:ModeKnn}
{\rm mode}_{N,N}(\ell):=\max_{1 \leq \ell \leq N^2} F_{N,N}(\ell) 
\approx 0.74118\,N^2.
\end{equation}
as depicted in Figure~\ref{fig:DistributionLengthNN}. We observe that the distribution is negatively skewed, the mean length of these paths being larger than the most frequent length,
\[\langle \ell\rangle_{N,N} :=\frac{\sum_{\ell=1}^{N^2} \ell\,F_{N,N}(\ell)}{T(2N+1,N+1)}
\approx 0.8125\,N^2 > {\rm mode}_{N,N}(\ell).
\] 
The above estimations were obtained by using a relatively low (N=8) dimension. As mentioned above, already at this low dimension we obtain the accurate qualitative behavior of the asymptotic distribution. In this way we can qualitatively describe a typical synchronization path for the Laplacian of the complete bipartite graph, starting at a random balanced ordered initial condition of dimension $2N$. For instance, such a synchronization path would most likely be of the length indicated in Equation~(\ref{eq:ModeKnn}).

\section{Remark and comments}~\label{sec:Final}

\ms Thanks to the monotonic behavior of the Laplacian flow in $K_N$, it was possible to completely describe the behavior of the transient dynamics of the system via a codification of the synchronized subnetworks by increasing functions. On the other hand, in the case of $K_{N,N}$, a similar codification is limited only to synchronizing paths starting at balanced initial conditions, which are the ones for which a monotonous behavior is obtained.\\

In both cases we obtained a closed formula for the number of number of realizable states, states given by combinatorial objects codifying all the realizable synchronized subnetworks. Moreover, the number of paths towards the synchronization of the two systems, which can be seen as a complexity function, remains an open problem. We can nevertheless obtain bounds that give us an idea of their growth order.\\

The probability density functions of the asymptotic distribution of the normalized length of a path towards synchronization in both cases are continuous, unimodal, and negatively skewed. The typical length with respect to the longest path is larger for $K_{N,N}$ than for $K_N$.\\

\ms Although the above results concern the Laplacian flow, they apply in some extend to the Kuramoto flow. In particular, in the case of complete network $K_N$, the transition diagram obtained from the Laplacian flow describes most of the paths towards synchronizations stating in a neighborhood of the diagonal. Inside this neighborhood we can use the coding of $\epsilon$-synchronized subnetworks defined for the Laplacian flow in Section~\ref{sec:KN} since the order of the coordinates is preserved by the Kuramoto flow, and therefore the increasing functions in $\Phi_N$ are suitable for the coding. Indeed, according to Equation~\eqref{eq:KuramotoModel} we have
\begin{eqnarray*}\label{eq:MonotonKuramoto}
\frac{d(x_n-x_m)}{dt}
&=&
\sigma\left(\sum_{k=1}^N\sin(x_k-x_n)-\sin(x_k-x_m)\right),\\
&=&\sigma\,R\left(\sin(\Theta-x_n)-\sin(\Theta-x_m)\right),
\end{eqnarray*}
where $R\,e^{i\Theta}=\left(\sum_{k=1}^N\cos(x_k)\right)+i\left(\sum_{k=1}^N\sin(x_k)\right)$. Hence, whenever $x_n=x_m$, $d(x_n-x_m)/dt=0$, which implies that the order in the coordinates is preserved under the flow since no crossing of coordinates is possible. Let us assume that $\max\{|x_n-\bar{x}|:\,1\leq n\leq N\}<\pi/4$, where $\bar{x}=\sum_{n=1}^Nx_n(0)$. In this case $|\Theta-\bar{x}|\leq \pi/4$ and $d(x_n-x_m)/dt=0$ if and only if $x_m=x_n$. Furthermore, the sign of $\sin(\Theta-x_n)-\sin(\Theta-x_m)$ is in this case the same as the sign of $x_m-x_n$, and therefore $|x_n-x_m|$ decreases monotonously for all initial condition. We have performed some numerical experiments and verify that the transition diagram defined in Section~\ref{sec:KN} is respected by the Kuramoto flow if one considers $\epsilon$ sufficiently small with respect to $\pi/4$ and initial conditions $x\in (S^1)^V$ such that $|x_n-\bar{x}|<\pi/4$  for all $1\leq n\leq N$. 

For $K_{N,N}$, the order of the coordinates at each of the two parties is preserved by the Kuramoto flow. For this we proceed as in the previous paragraph and obtain 
\begin{eqnarray*}\label{eq:MonotonKuramotoNN}
\frac{d(x_n-x_m)}{dt}&=&\sigma\,R_2\left(\sin(\Theta_2-x_n)-\sin(\Theta_2-x_m)\right),\\
\frac{d(x_{N+n}-x_{N+m})}{dt}&=&\sigma\,R_1\left(\sin(\Theta_1-x_n)-\sin(\Theta_1-x_m)\right),
\end{eqnarray*}
where $R_1\,e^{i\Theta_1}=\left(\sum_{k=1}^N\cos(x_k)\right)+i\left(\sum_{k=1}^N\sin(x_k)\right)$ and similarity for $R_2\,e^{i\Theta_2}$. From this it follows that if $x_n=x_m$ then $d(x_n-x_m)/dt=0$, and similarly for $x_{N+n}-x_{N+m}$. Therefore the order in the coordinates at each party is preserved under the flow which allows us to use the coding of $\epsilon$-synchronized subnetworks defined for the Laplacian flow in Section~\ref{sec:KNN}. 

\ms As already mentioned, the transition diagram defined for the Laplacian flow over $K_{N,N}$ describes only the paths toward synchronization corresponding to balanced initial conditions. In Figure~\ref{fig:N2TransitionGraph} we marked in red the subnetworks incompatible with balanced initial conditions. The whole transition diagram, which contains those subnetworks, admits non-monotonous paths. Furthermore, for unbalanced initial conditions, the order in the differences between coordinates is not preserved by the flow. The description of the full transition diagram for the Laplacian flow over $K_{N,N}$, would be the subject of future work.\\

\ms Finally we would like to emphasize that these synchronizing sequences can be seen as partitioning the basin of attraction of a given attractor (here the fully synchronized state). Since for a given finite $\epsilon$ the final synchronized network will be reached in a finite time $\tau(\epsilon,N)$, if  the space of initial conditions has a finite volume, the full space-time will be as well bounded, and these sequences are partitioning that full space time. Moreover, by associating to a given sequence  an ensemble of initial conditions realizing that sequence, we should be able to  measure that ensemble and add corresponding weights (measures) to each sequence and characterize even further the space-time complexity.

\section*{Acknowledgements}
\noindent The authors are grateful for the financial support from ECOS-CONACyT-ANUIES through Grant M16M01. A. E. benefits from the National Scholarship No. 722957 offered by CONACyT-M\'exico. She thanks Fundaci\'on Sofia Kovalevskaia and Sociedad Matem\'atica Mexicana for financial support. E. U. benefited from a research stay at the Centre de Physique Th\'eorique-Luminy, financed by the CNRS-France.

\bibliography{mybibfile}

\appendix

\section{}\label{app:KNIncreasing}

\ms For each increasing function $\phi:\{1,2,\ldots,N\}\to\{1,2,\ldots,N\}$ such that $\phi\geq {\rm Id}$, there exists an ordered initial condition $x\in\IR^N$ such that $\phi=\phi_x$. For this we use a representation of $\phi$ as a disjoint union of directed trees as follows. Let ${\rm Fix}(\phi):=\{1\leq n\leq N:\, \phi(n)=n\}$. To each $n\in{\rm Fix}(\phi)$ we associate a directed tree $T_n$, rooted at $n$, with vertex set $V_n:=\bigcup_{l=0}^{h(n)}\phi^{-l}(\{n\})$ and directed edges in $A_n:=\{(k,\phi(k)):\, k\in V_n\setminus\{n\}\}$. The vertex set $V_n$ splits into $h(n)+1$ disjoint levels, $V_n^l:=\phi^{-l}(\{n\})$, $0\leq l\leq h(n)$. The number $h(n)$ is the high $T_n$. The maximal paths in $T_n$ are completely determined by their starting vertices, which have to be leaves. Let $\ell_n^1 < \ell_n^2 < \cdots < \ell_n^{w(n)}$ be the leaves of $T_n$. Its number, $w(n)$, is the width of the tree $T_n$. Since $\phi$ is increasing and such that $\phi \geq {\rm Id}$, then every element in the $l$-th level, $V_n^{l}$, is greater than all the elements in the $l'$-th level, $V_n^{l'}$ whenever $l < l'$. It implies that the length $l(m)$ of the path starting at $m$ and ending at the root, is a decreasing function of $m$. Each maximal path in $T_n$ starts at a leaf and the longest of those paths have length $h(n)$, and start at leaves in the highest level. Furthermore, all vertices in $T_n$ belong to a maximal path, which means that it is reachable from a leaf.

\ms Now, given $\phi:\{1,2,\ldots,N\}\to\{1,2,\ldots,N\}$ increasing and such that $\phi\geq{\rm Id}$, let $\{T_{n_k}:\, 1\leq k\leq R\}$ be the associated collection of directed trees and $n_1 <n_2 < \cdots < n_R$ in ${\rm Fix}(\phi)$ the corresponding roots. Define $x\in \IR^N$ such that $x_{n_1}=\epsilon\,h(n_1)$, and for each $1 \leq k < R$, 
\begin{equation}\label{eq:SeparationTrees}
x_{n_{k+1}}=x_{n_k}+(h(n_k)+2)\,\epsilon.
\end{equation} 
In this way, we fix the value of $x_n$ for each $n\in {\rm Fix}(\phi)$ in such a way that 
$x_{n_k}+\epsilon < x_{n_{k+1}}-h(n_{k+1})\,\epsilon$ for each  $1\leq k < R$. Now, for each $n\in {\rm  Fix}(\phi)$, let $\ell_n^1 < \ell_n^2 < \cdots < \ell_n^{w(n)}$ be the leaves of $T_n$. For each $1 \leq j \leq w(n)$ and $0\leq k\leq l(n_j)$ for which $x_{\phi^k(\ell_n^j)}$ is not yet defined, let 
\begin{equation}\label{eq:InsideTrees}
x_{\phi^k(\ell_n^j)}=x_n-(l(n_j)-k)\,\epsilon+(j-1)\frac{\epsilon}{w(n)}.
\end{equation} 
Let us remind that $l(n_j)$ is the length of the maximal path starting at $\ell_n^j$. It is not difficult to verify that Equations~\eqref{eq:SeparationTrees} and~\eqref{eq:InsideTrees} define an ordered initial condition $0=x_1<x_2<\cdots<x_N=\sum_{k=1}^R(h(n_k)+2)$, such that $\phi_x=\phi$. 

\section{}\label{app:KNNIncreasing}

\ms Each couple of increasing functions $\alpha,\omega:\{1,\ldots,N\}\to\{0,1,\ldots,N+1\}$ is compatible with some $x\in \IR^{2N}$ in terms of the Equations~(\ref{eq:KNNIncreasingFunctionA}) and~(\ref{eq:KNNIncreasingFunctionB}), and therefore codify an $\epsilon$-synchronized subnetwork, provided ${\rm im}(\alpha)\subset[1,N+1]$, ${\rm im}(\omega)\subset[0,N]$ and $\alpha\leq \omega+1$. Such an initial condition can be constructed as follows.

\ms For each $1\leq n\leq N$ let $\A_n:=\{1\leq m\leq N:\, \alpha(n)\leq m\leq\omega(n)\}$. Let us partition $\{1,2,\ldots,N\}=\bigsqcup_{k=1}^\ell I_k$, where for each $1\leq k\leq\ell$, $I_k=\{n_k,n_k+1,\ldots,m_k\}$ is such that $\A_n\cap\A_{n+1}\neq \emptyset$ for each $n_k\leq n < m_k$ and it is a maximal element in the sense of inclusion ($I_k \subsetneq I\Rightarrow \bigcup_{n\in I_k}\A_n$ is not an interval). Notice that $n_1=1$ and that $I_k=\{n_k\}$ whenever $\alpha(n_k)=\omega(n_k)+1$. 

\ms For each $1 < k \leq \ell$, let $\Delta:I_k\to I_k$ be such that $\Delta(n)=\max\{m\in I_k:\, \A_n\cap\A_m\neq\emptyset\}$. Clearly $\Delta(n)\geq n$ and $\Delta(n)=n$ if and only if $n=n_k=m_k$. We can associate to $\Delta$ a directed tree $T_k$ with vertices in $I_k$, rooted at $m_k$, and arrows $n\mapsto\Delta(n)$. The structure of these trees is similar to that of the trees described in Appendix~\ref{app:KNNIncreasing}. Let $n_k\mapsto \Delta(n_k)\mapsto\cdots\mapsto \Delta^j(n_k)\mapsto\cdots\mapsto m_k =\Delta^{h_k}(n_k)$ be the maximal path in $T_k$ and for each $1\leq j\leq l_k$ let $V_j=\Delta^{-j}(\{m_k\})$ be the $j$-th level of $T_k$. Clearly $\min V_j=\Delta^{h_k-j}(n_k)$ and $\max V_j < \min V_{j-1}$ for each $0\leq j\leq h_k$. 

\ms Assume $x_{n_k}$ is given. Let $n_{k,j}:=\min V_j$ and define $x_{n_{k,j}}:=x_{n_k}+j\epsilon$ for each $1\leq j\leq h_k$.  Now, for $n_{k,j}\leq n < n_{k,j-1}$, let $x_n=x_{n_{k,j}}+(n-n_{k,j})\,\epsilon/(n_{j-1}-n_{k,j})$. With $\delta_k:=\frac{1}{2}\min_{n_k\leq n < m_k}(x_{n+1}-x_n)$, for each $n_k\leq n < m_k$ and $\alpha(n)\leq m <\alpha(n+1)$, let $x_{N+m}=x_n-(\epsilon-\delta_k)$. For $n_{k,1}\leq n < n_{k,0}\equiv m_k$ and $\omega(n) < m \leq\omega(n+1)$, let  $x_{N+m}=x_n+(\epsilon-\delta_k)$. Finally, for $\alpha(m_k)\leq m \leq \omega(n_{k,1})$, define $x_{N+m}=(x_{n_{k,1}}+x_{m_k})/2$.

\ms In order to complete the specification of all the coordinate, fix $x_1 = x_{n_1}=0$ and for each $1\leq k\leq \ell$ let $x_{n_k}:=x_{m_{k-1}}+3\epsilon$. Finally, for each $m \notin \bigcup_{n=1}^N\A_n$, let $k(m):=\min\{1\leq k\leq \ell:\, \alpha(n_k) > m\}$ and define $x_{N+m}:=x_{N+\alpha(n_k)}-3\epsilon/2$. If $\omega(N)<N$, then define $x_{N+m}:=x_{m_\ell}+3\epsilon/3$.

\end{document}